\documentclass{amsart}
\usepackage{hyperref,amsmath,amsfonts,amssymb,amsthm,enumerate,mathtools,graphicx,xparse,xcolor,tikz,scalerel,stackengine}
\usepackage[a4paper, left=2cm, right=2cm, top=3cm, bottom=2cm]{geometry}
\usepackage[utf8]{inputenc}
\usepackage{enumitem}
\usepackage{float}
\usepackage{mdframed}

\graphicspath{{figures/}}
\usepackage{pgfplots}
\pgfplotsset{compat=1.18}
\usepackage{subfig}
\usepackage{esint}    

\usepackage{accents}
\usepackage{empheq}
\usepackage{subfig}
\usepackage{multimedia}
\usepackage{animate}
\usepackage[makeroom]{cancel}
\definecolor{morange}{rgb}{0.92, 0.51, 0.11}
\definecolor{mLightGreen}{HTML}{14B03D}
\definecolor{mPink}{HTML}{e1388d}
\definecolor{mRed}{HTML}{E83C3C}

\usepackage{scalerel}

\numberwithin{dummy}{section}
\numberwithin{equation}{section}

\DeclareMathOperator*{\essup}{ess\,sup}

\begin{document}

\theoremstyle{plain}
\newtheorem{thm}{Theorem}[section]
\newtheorem{lemma}[thm]{Lemma}
\newtheorem{claim}{Claim}
\newtheorem*{cor}{Corollary}
\newtheorem{prop}{Proposition}

\theoremstyle{plain}
\newtheorem{df}[thm]{Definition}

\theoremstyle{remark}
\newtheorem*{rem}{Remark}
\newtheorem*{example}{Example}
\newtheorem*{conv}{Convention}


\newenvironment{myproof}{
  \par\medskip\noindent
  \textit{Proof}.
}{
\newline
\rightline{$\qedsymbol$}
}


\newcommand{\rr}{\mathbb{R}}
\newcommand{\nn}{\mathbb{N}}
\newcommand{\om}{(\Omega)}
\newcommand{\V}{W^{1,p}_{0,\diver}}
\newcommand{\Leb}{L^{2}_{0,\diver}}
\newcommand{\boch}[2]{#1(0,T;#2)}
\newcommand{\bdary}{\partial\Omega}

\newcommand{\summing}[2]{\sum_{#1 = 1}^#2}
\newcommand{\intom}{\int_\Omega}
\newcommand{\intq}{\int_Q}
\newcommand{\intt}{\int_{0}^T}
\newcommand{\inttas}{\int_{0}^{T^\ast}}
\newcommand{\intu}[1]{\int_{0}^#1}
\newcommand{\intd}[1]{\int_{#1}}
\newcommand{\dx}{\,\mathrm{d}x}
\newcommand{\dt}{\,\mathrm{d}t}
\newcommand{\dtau}{\,\mathrm{d}\tau}
\newcommand{\dq}{\, \mathrm{d}x\, \mathrm{d}t}
\newcommand{\dtom}{ \dx \dtau}
\newcommand{\ds}{\,\mathrm{d}s}

\newcommand{\goto}{\rightarrow}
\newcommand{\gotow}{\rightharpoonup}
\newcommand{\gotows}{\rightharpoonup^\ast}

\newcommand{\abs}[1]{\left|{#1}\right|}

\newcommand{\diver}{\operatorname{div}}
\newcommand{\ddt}{\frac{d}{dt}}
\newcommand{\dert}{\partial_t}
\newcommand{\Tr}{\textrm{Tr }}

\newcommand{\parcd}[2]{\frac{\partial #1}{\partial #2}}
\newcommand{\dual}[2]{\left\langle#1, \,#2\right\rangle}
\newcommand{\iprodS}[3]{\langle #1, #2 \rangle_{#3}}
\newcommand{\norm}[3]{\|{#1}\|_{#2}^{#3}}
\newcommand{\normg}[3]{\|{\nabla #1}\|_{#2}^{#3}}
\newcommand{\norms}[3]{\|{#1}\|_{1,#2}^{#3}}

\newcommand{\vth}{\vartheta}
\newcommand{\htet}{\hat{\vartheta}}
\newcommand{\bth}{\pmb \vartheta}
\newcommand{\bphi}{\pmb \varphi}

\newcommand{\bU}{{\pmb U}}
\newcommand{\bu}{\pmb u}
\newcommand{\bv}{\pmb v}
\newcommand{\bw}{\pmb w}
\newcommand{\bq}{\pmb q}
\newcommand{\bp}{\pmb p}
\newcommand{\beh}{\pmb f}

\newcommand{\Str}{\boldsymbol{\mathcal{T}}}
\newcommand{\tens}{\boldsymbol{\mathcal{S}}}
\newcommand{\tk}{\mathcal{T}_k}
\newcommand{\tke}{\mathcal{T}_{k,\varepsilon}}
\newcommand{\gk}{\mathcal{G}_k}
\newcommand{\dens}{\boldsymbol{\mathcal{D}}}
\newcommand{\bO}{\boldsymbol{{0}}}
\newcommand{\wens}{\boldsymbol{\mathcal{W}}}

\newcommand{\upnm}{^{N,M}}
\newcommand{\upnn}{^N}
\newcommand{\downdiv}{_{0,\diver}}
\newcommand{\upsob}{^{1,p}}

\newcommand{\todo}[1]{\color{red}TODO: #1 \color{black}}
\newcommand{\quest}[1]{\color{blue}QUESTIONABLE: #1 \color{black}}

\newcommand{\cQ}{{\mathcal Q}}
\newcommand{\cH}{{\mathcal H}}
\newcommand{\cT}{{\mathcal T}}
\newcommand{\dd}{\operatorname{d\!}}

\title[Time continuity of solutions to Navier--Stokes--Fourier system] {On the continuity in time of solutions to a generalized Navier--Stokes--Fourier system}\thanks{This work has been supported by the project 25-16592S financed by Czech science foundation (GA\v{C}R) and by Charles University Centre program No. UNCE/24/SCI/005. M.~Bul\'{\i}\v{c}ek and P.~Kaplick\'{y} are  members of the Ne\v{c}as Center for Mathematical Modeling. }
\author[M.~Bul\'{i}\v{c}ek]{Miroslav Bul\'{\i}\v{c}ek}
\address{Charles University, Faculty of Mathematics and Physics, Mathematical Institute\\
Sokolovsk\'{a} 83, 186 75 Prague, Czech Republic}
\email{mbul8060@karlin.mff.cuni.cz}

\author[P. Kaplick\'{y}]{Petr Kaplick\'{y}}
\address{Charles University, Faculty of Mathematics and Physics, Department of Mathematical Analysis \\
Sokolovsk\'{a} 83, 186 75 Prague, Czech Republic}
\email{kaplicky@karlin.mff.cuni.cz}

\author[L. Wintrov\'{a}]{Lucie Wintrov\'{a}}
\address{Charles University, Faculty of Mathematics and Physics, Mathematical Institute \\
Sokolovsk\'{a} 83, 186 75 Prague, Czech Republic}
\email{wintrova@karlin.mff.cuni.cz}

\begin{abstract}
We consider the flow of a generalized non-Newtonian incompressible heat-conducting fluid in a~bounded two-dimensional domain, subject to Dirichlet boundary conditions for velocity and temperature. The fluid obeys a power-law constitutive relation for the Cauchy stress with exponent~$p$. For $p\geq 2$ and finite-energy initial data, we establish the existence of a global-in-time weak solution that satisfies the entropy equality.
The novelty of this work is the rigorous proof of time continuity of the temperature in $L^1(\Omega)$, a property not previously established in this setting. 

Furthermore, we prove regularity and time continuity for a weak solution of the entropy equation with a convective term and an $L^1$ right-hand side under minimal assumptions on the velocity regularity, in arbitrary spatial dimensions. We show that this continuity is equivalently described by vanishing dissipation on high level sets, a truncated variational inequality for admissible test functions, or the associated equality. This reveals the connection between energy dissipation, weak stability, and temporal regularity.

\end{abstract}

\maketitle




\section{Introduction}\label{sec:1}
In this paper, we consider the following system of partial differential equations
\begin{align} \label{i1}
\dert \bu + \diver (\bu \otimes \bu)
- \diver \tens + \nabla \pi
	&= \beh,
\\ \label{i2}
\diver \bu &= 0,
\\ \label{i3}
\dert \vth + \diver (\vth \bu)
+ \diver \bq 
	&= \tens : D\bu,
\\ \label{i4}
\dert \eta + \diver (\eta \bu)
+ \diver \left(\frac{\bq}{\vth}\right) 
	&= \frac{\tens : D\bu}{\vth} - \frac{\bq \cdot \nabla \vth}{\vth^2},
\end{align}
which is satisfied in $Q := (0,T) \times \Omega$, where \(T>0\) and \(\Omega \subset \mathbb{R}^2\) is a bounded Lipschitz domain.
Here $\bu: Q\to \rr^2$ denotes the velocity field, 
$D\bu:= (\nabla \bu + (\nabla \bu )^T) /2 $ is the symmetric part of the velocity gradient,  $\pi:Q\to \rr$ is the pressure,  $\vth: Q\to \rr$ the temperature, and $\eta = \log \vth$ the entropy. Furthermore, $\beh:Q \to \rr^2$  is a given density of the external body forces, $\tens: Q \to \rr^{2\times 2}$ denotes the viscous part of the Cauchy stress tensor, and $\bq: Q \to \rr^2$ is the heat flux. The system~\eqref{i1}--\eqref{i4} is completed by the initial  and
boundary conditions of the form
\begin{equation}
    \begin{aligned} \label{ic}
    \bu &= \pmb{0}, 
&\vth = \vth_{b} &\;\;\;\textrm{on $\partial \Omega\times (0,T)$},
\\
\bu(0) &= \bu_0,
&\vth(0) = \vth_0 &\;\;\;\textrm{in $\Omega$}.
\end{aligned}
\end{equation}
The above system describes the planar flow of an incompressible, homogeneous, heat-conducting fluid, where a unit density is assumed for simplicity. Equation \eqref{i2} represents the incompressibility constraint. The balance of linear momentum is given by \eqref{i1}, and the balance of internal energy by \eqref{i3}, where we assume that the heat capacity is identically equal to one. Equation \eqref{i4} describes the balance of entropy. At first glance, the system might appear to be overdetermined; however, this is not the case, since \eqref{i3} and \eqref{i4} are formally equivalent, provided that the solution is sufficiently smooth.

We impose homogeneous Dirichlet boundary conditions for the velocity and spatially inhomogeneous Dirichlet boundary conditions for the temperature. Consequently, no mass passes through the boundary, although thermal interaction with the exterior is allowed.

\subsection{Constitutive relations}\label{constrel}

We allow the material parameters to depend on the temperature and assume that the fluid exhibits non-Newtonian behavior. Accordingly, the viscous part of the Cauchy stress tensor $\tens$ is taken to depend nonlinearly on the symmetric part of the velocity gradient. This dependence is modeled by a power law with index $p \geq 2$.  More precisely, we assume that $\tens = \tens^*(\vth, D\bu)$, where $\tens^*: (0, \infty) \times \rr_{sym}^{2\times2} \to \rr_{sym}^{2\times2}$ is a continuous mapping, and that there exist constants $0 < \underline{\nu} \le \overline{\nu} < \infty$ and $p \ge 2$ such that, for all $\vth \in \rr_+$, $\dens_1, \dens_2 \in \rr_{sym}^{2\times2}$, it holds
\begin{equation}\label{tensor}
    \begin{aligned}
        (\tens^*(\vth, \dens_1) - \tens^*(\vth, \dens_2)):(\dens_1 - \dens_2) &\geq 0,\\
        \tens^*(\vth, \dens_1):\dens_1 &\geq \underline{\nu}|\dens_1|^p-\overline{\nu},\\
        |\tens^*(\vth, \dens_1)| &\leq \overline{\nu}(1+|\dens_1|)^{p-1} \quad \textrm{ and } \quad \tens^*(\vth, \bO)= \bO.
    \end{aligned} 
\end{equation}
The above assumptions are sufficient in cases where we deal with the existence and continuity of the solution. In the case where we also want to show its exponential stability, we additionally require that
\begin{equation}\label{stability}
        \tens^*(\vth, \dens_1): \dens_1 \ge \underline{\mu}(1+|\dens_1|)^{p-2}|\dens_1|^2.
\end{equation}

Next, we specify the constitutive relation for the heat flux, which is assumed to obey Fourier’s law. More precisely, we set
$\bq = \bq^\ast(\vth)$ where
\begin{equation}\label{hf}
    \bq^\ast(\vth) = -\kappa(\vth)\nabla\vth,
\end{equation}
and the heat conductivity $\kappa: \rr \to (0, \infty)$ is a continuous function. We further assume that there exist constants $0<\underline{\kappa}$, $\overline{\kappa}<\infty$ such that
\begin{equation}\label{kap}
   \underline{\kappa}\leq \kappa \leq \overline{\kappa}. 
\end{equation}


\subsection{Notation}\label{not}
In what follows, we adopt the standard notation for the Lebesgue, Sobolev, and Bochner spaces, equipping them with their usual norms. The symbol $C^\infty_0$ is reserved for smooth, compactly supported functions. The function space associated with the incompressible setting is denoted by 
$$
\V:= \{\bv \in W\upsob_0 (\Omega;\rr^2); \diver \bv = 0\},
$$ 
and $\Leb$ denotes its closure in the $L^2$ topology. For a general Banach space, the duality pairing is denoted by $\dual{\cdot}{\cdot}_X$, while in the case $X=\V$ the subscript is omitted for simplicity. For notational convenience, we also use the abbreviation $V:=\V$.

We recall here a few classical inequalities related to the function space setting used in this paper. First, we recall the Korn inequality  
\begin{equation}\label{korn}
      C\norm{\bu}{1,q}{ } \leq \norm{D\bu}{q}{ } + \norm{\Tr \bu}{L^2(\partial \Omega)}{ },
\end{equation}  
which holds true for any bounded Lipschitz domain $\Omega\subset \rr^2$ and $q \in (1,\infty)$, see \cite[as Lemma 1.11]{BuMaRa07}. Note that for the space $V$, the trace part of the above inequality vanishes.

Second is the classical Lebesgue--Sobolev interpolation inequality 
\begin{equation}\label{galnir}
       \|z\|_r \leq C\norm{z}{1,s}{\frac{2s(q-r)}{r(2q-sq-2s)}}\norm{z}{q}{\frac{q(2r-sr-2s)}{r(2q-sq-2s)}},
\end{equation}
which holds true for $r \ge q \ge 1$ as follows. If $s < 2$, then the inequality is valid for all $r \leq \frac{2s}{2-s}$. If $s > 2$, then it holds for all $r \leq \infty$. If $s = 2$, then it holds for all $r < \infty$.

\subsection{Assumptions on data}
Here we state the assumptions imposed on the data. We begin with $\vth_b$, for which we require
\[
\vth_b \in W^{\frac12,2}(\partial \Omega) \cap L^{\infty}(\partial \Omega).
\] 
Note that, by classical elliptic theory and \eqref{hf}--\eqref{kap}, this assumption guarantees the existence of a function $\hat{\vth}: \Omega \to \rr$ satisfying
\begin{align}\label{hatcond}
     &\diver\bq^\ast(\htet) = 0 \quad \text{in $\Omega$,} && \htet = \vth_b \quad \textrm{on $ \bdary$,} && \hat{\vth} \in L^\infty(\Omega) \cap W^{1,2}(\Omega).
\end{align}
Hence, in what follows, we impose assumptions only on $\htet$ and no longer consider $\vth_b$.  

We consider initial and external data satisfying standard assumptions. In particular, the external force $\beh$ and the initial velocity $\bu_0$ are assumed to fulfill
\begin{equation}
\label{conditions}
        \beh \in L^{p'}\left(0,T;V^\ast\right), \quad \bu_0 \in \Leb\om,
\end{equation}
while the initial and boundary temperatures are assumed to satisfy
\begin{gather} 
        \label{conditionstheta} \vth_0 \in L^1\om, \quad \htet \in W^{1,2}\om \cap L^\infty\om, \\ 
        \label{muconditions}
        \mu := \min \left\{\mathop{\operatorname{ess\,inf}}\limits_{x \in \Omega}  \hat{\vth}(x), \mathop{\operatorname{ess\,inf}}\limits_{x \in \Omega} \vth_0(x) \right\} > 0.
\end{gather}

\subsection{Main results}\label{formulation}

The main results of the paper are summarized in the following theorems.

\begin{thm}[On the existence of a solution fulfilling entropy equality]\label{maintheorem}
Let $\Omega \subset \rr^2$ be a bounded domain with Lipschitz boundary and $T>0$. Assume that $\tens^\ast$ and $\kappa$ satisfy \eqref{kap}--\eqref{tensor}. Additionally, assume that $\beh$, $\bu_0$, $\vth_0$, and $\htet$ fulfill \eqref{conditions}--\eqref{muconditions}.
    
     Then there exists a quadruplet $(\bu, \tens, \vth, \eta)$ fulfilling
    \begin{align}
        \bu &\in C([0,T];L^2_{0, \diver}) \cap L^p(0,T;V),
        \\
        \dert \bu &\in  L^{p'}(0,T;V^\ast), \; \tens \in  L^{p'}(Q, \rr^{2\times2}_{sym}),
        \\
        \vth &\in C([0,T];L^1\om), \; (\vth)^\alpha \in L^2(0,T; W^{1,2}\om) && \textrm{ for all }\alpha \in \left[0,{1}/{2}\right),
        \\
        \vth &\in L^r(Q) && \textrm{ for all }r \in \left[1,2\right),
        \\
        \vth - \hat{\vth} &\in L^s(0,T; W^{1,s}_0\om) && \textrm{ for all }s \in \left[1,{4}/{3}\right),
        \\
        \eta &\in L^2(0,T; W^{1,2}\om)\cap L^q(Q) && \textrm{ for all }q \in [1,\infty),
    \end{align}
    and satisfying \eqref{i1}--\eqref{ic} in the following sense:
    \\
    \emph{Momentum equation:} The Cauchy stress is of the form $\tens = \tens^\ast(\vth, D\bu)$ a.e. in $Q$, it holds $\bu(0)=\bu_0$ in $L^2$, and for all $\bw \in L^p(0,T; V)$
\begin{equation}\label{me}
    \begin{aligned}
        \intt \dual{\dert \bu}{\bw} \dt &+  \intt \intom \tens:D\bw \dx \dt \\
    &= \intt\intom (\bu \otimes \bu):D\bw \dx \dt + \intt \dual{\beh}{\bw} \dt ;
    \end{aligned} 
\end{equation}
\\
\emph{Internal energy balance:} Temperature satisfies the minimum principle $\vth \geq \mu$ a.e. in $Q$, it holds $\vth(0) = \vth_0$ in $L^1\om$, and for all $\varphi \in C^\infty_0((-\infty,T) \times \Omega)$
\begin{equation}\label{ieb}
    \begin{aligned}
        -\int_{0}^T \intom \vth \dert \varphi \dx \,dt - \int_{0}^T \int_\Omega \vth \bu \cdot \nabla \varphi
        \dx \,dt + \int_{0}^T \intom \kappa(\vth) \nabla\vth \cdot \nabla \varphi \dx \,dt \\
    = \int_{0}^T\intom \tens:D\bu \varphi \,dx \,dt + \intom \vth_0 \varphi(0) \dx;
    \end{aligned} 
\end{equation}
\\
\emph{Entropy equality:} Entropy is given as $\eta = \ln \vth$ a.e. in $Q$, $\eta_0 := \ln \vth_0$ and for all $\varphi \in C^\infty_0((-\infty,T) \times \Omega)$
\begin{equation}\label{ee}
    \begin{aligned}
        -\intt \intom \eta \dert \varphi \dx \dt - \intt 
        \intom \eta \bu \cdot \nabla \varphi
        \dx \dt + \intt \intom \kappa(\vth) \nabla\eta \cdot \nabla \varphi \dx \dt \\
    = \intt \intom \frac{1}{\vth} \tens:D\bu \varphi \dx \dt + \intt \intom \kappa(\vth)
    \frac{|\nabla\vth|^2}{\vth^2} \varphi \dx \dt + \intom \eta_0 \varphi(0) \dx.
    \end{aligned} 
\end{equation}
\end{thm}

\bigskip

In the second main result, we discuss properties of weak solutions of the entropy equation with a convective term and temperature-dependent conductivity, where the velocity field and the right-hand side have only minimal regularity. The result holds in arbitrary spatial dimensions and builds on the ideas developed for the first main theorem. We show that the entropy equality is a key property and that it is equivalent to the renormalized formulation of the temperature equation. In addition, we provide a necessary and sufficient condition under which $\vth \in C([0,T]; L^1(\Omega))$. While Theorem \ref{maintheorem} establishes this result in the two-dimensional setting for a particular right-hand side, we now generalize it to arbitrary dimensions and to any right-hand side in $L^1(Q)$.

Before stating the second main theorem, we introduce the required notation. For any $m>0$ we define
\begin{align}\label{tk}
    &\mathcal{T}_m:\rr \to [-m,m], &\mathcal{T}_m(z) := \text{sign} (z) \min\{|z|,m\}.
\end{align}
For $\delta \in (0, m)$, let $\mathcal{T}_{m,\delta} \in C^2(\rr)$ denote the mollification of $\mathcal{T}_m$, which is given by convolution with $\omega_\delta(\cdot)=\omega(\cdot/\delta)/\delta$ where $\omega$ is a smooth symmetric regularization kernel supported in $[-1,1]$ such that $\omega$ is positive in $(0,1)$. For $\mathcal{T}_{m,\delta}$ we have
\begin{equation}\label{propertiestke}
    \begin{aligned}
    &\mathcal{T}_{m,\delta}(z)=\mathcal{T}_{m}(z) &&\text{if } |z| \leq m-\delta \text{ or } |z| \geq m+\delta,
    \\
    & 0 \leq \mathcal{T}_{m,\delta}' \leq 1, \quad   \mathcal{T}_{m,\delta} \leq \mathcal{T}_{m} &&\text{on }(0,\infty),\\
    & -\frac{C}{\delta}\leq
    \mathcal{T}_{m,\delta}''=-\frac1\delta\omega(\frac{\cdot-m}\delta)\leq -\frac c\delta\chi_{(m-\frac\delta2,m+\frac\delta2)}\leq 0,  \quad
    &&\text{on $(0,+\infty)$}.
\end{aligned}
\end{equation}
Further, we define the primitive function 
\begin{equation}\label{dfGk}
\mathcal{G}_k(t) := \int_0^t \mathcal{T}_k(\tau)\, d\tau 
= \frac{t^2}{2}\chi_{\{|t|\le k\}} +  \left(k|t| - \frac{k^2}{2}\right) \chi_{\{|t|>k\}}.
\end{equation}

Finally, we can formulate the second main result of the paper.
\begin{thm}[On the role of entropy equality]\label{maintheoremB}
Let $\Omega \subset \rr^d$ be a bounded domain with Lipschitz boundary and $T>0$. Assume that $\kappa$ satisfies \eqref{kap}, $\htet$ and $\vth_0$ satisfy \eqref{conditionstheta}--\eqref{muconditions} and that $g\in L^1(Q)$. 
Let
\begin{align}
\label{vb2B}
\vth &\in L^{\infty}(0,T;L^1\om),\\
\bu &\in L^p(0,T;L^p_{0, \diver}) \label{vb1}
\end{align}
for some $p>2$, and assume that there exists a constant ${\mu}>0$ such that 
$$
\vth \ge \mu.
$$
Let us define the entropy $\eta:=\ln \vth$, the representation of the boundary value $\hat{\eta}:=\ln \htet$ and its initial value $\eta_0:= \ln \vth_0$, and assume that
\begin{align}
\eta-\hat{\eta} &\in L^2(0,T; W^{1,2}_0\om)\cap L^q(Q) \qquad \textrm{for all } q \in [1,\infty), \label{vb4}
\end{align}
and that for all $\varphi \in C^1_0((-\infty,T)\times \Omega)$ the entropy equality is satisfied in the following sense
\begin{equation}
\label{eeB}
    \begin{aligned}
        -\intt \intom \eta \dert \varphi \dx \dt - \intt 
        \intom \eta \bu \cdot \nabla \varphi
        \dx \dt + \intt \intom \kappa(\vth) \nabla\eta \cdot \nabla \varphi \dx \dt \\
    = \intt \intom \frac{g\,  \varphi }{\vth} \dx \dt + \intt \intom \kappa(\vth)
    |\nabla \eta|^2 \varphi \dx \dt + \intom \eta_0 \varphi(0) \dx.
    \end{aligned} 
\end{equation}
Then we have the following regularity results 
    \begin{align}
        (\vth)^\alpha &\in L^2(0,T; W^{1,2}\om) &&\textrm{ for all } \alpha \in \left[0,\frac{1}{2}\right),\label{vb2}\\
        \vth &\in L^r(Q) &&\textrm{ for all } r \in \left[1,\frac{d+2}{d}\right),\label{vb3}\\
        \vth &\in L^s(0,T; W^{1,s}\om),&&\textrm{ for all } s \in \left[1,\frac{d+2}{d+1}\right),\label{tls}\\
        \vth^{\alpha} & \in C([0,T]; L^1(\Omega)) &&\textrm{ for all } \alpha \in [0,1), \label{contir}
    \end{align}
and we also have that $\vth(0)=\vth_0$, in the sense of the continuity result in \eqref{contir}. Moreover, the temperature satisfies the renormalized temperature equation, i.e., for all $k,\delta>0$ and all $\varphi \in C^1_0((-\infty,T)\times \Omega)$, it holds
\begin{equation}\label{iebR}
    \begin{aligned}
        -\int_{0}^T \intom \mathcal{T}_{k,\delta}(\vth) \dert \varphi \dx \dt - \int_{0}^T \int_\Omega \mathcal{T}_{k,\delta}(\vth) \bu \cdot \nabla \varphi
        \dx \,dt + \int_{0}^T \intom \kappa(\vth) \nabla \mathcal{T}_{k,\delta}(\vth) \cdot \nabla \varphi \dx \,dt \\
    = \int_{0}^T\intom g\, \mathcal{T}_{k,\delta}' (\vth) \varphi -\mathcal{T}_{k,\delta}''(\vth)\kappa(\vth)|\nabla \vth|^2 \varphi \dx \dt + \intom \mathcal{T}_{k,\delta}(\vth_0) \varphi(0) \dx.
    \end{aligned} 
\end{equation}
Furthermore, in the case that $p>(d+2)/2$, for all nonnegative $\varphi \in C_0^1((-\infty,T)\times \Omega)$, the temperature satisfies the following inequality
\begin{equation}\label{iebB}
    \begin{aligned}
        -\int_{0}^T \intom \vth \dert \varphi \dx \dt - \int_{0}^T \int_\Omega \vth \bu \cdot \nabla \varphi
        \dx \dt + \int_{0}^T \intom \kappa(\vth) \nabla\vth \cdot \nabla \varphi \dx \dt \\
    \ge \int_{0}^T\intom g \varphi \dx \dt + \intom \vth_0 \varphi(0) \dx.
    \end{aligned} 
\end{equation}
 The following statements are equivalent:
\begin{itemize}
\item[{\bf A)}]    
The temperature satisfies
\begin{equation} \label{Equi1}
\lim_{k\to \infty} \int_0^T \int_{\Omega} \frac{|\nabla \vth|^2}{\vth} \chi_{\{k\le \vth \le 2k \}} \dx \dt =0.
\end{equation}
\item[{\bf B)}] For every test function $\varphi \in L^2(0,T; W^{1,2}\om) \cap L^{\infty}(Q) \cap W^{1,1}(0,T; L^1\om)$ satisfying $\varphi-\vth\in L^1(0,T; W^{1,1}_0)$, and for all $k\in \mathbb{N}$ and a.e. $\tau\in (0,T)$, it holds
\begin{equation}
\label{Equi2}
\begin{aligned}
& \intom \mathcal{G}_k (\vth(\tau) -\varphi(\tau)) \dx +\int_{0}^{\tau}\intom \partial_t \varphi\, \mathcal{T}_k(\vth-\varphi) \dx \dt + \int_{0}^{\tau} \int_\Omega \nabla \varphi \cdot  \bu \mathcal{T}_k(\vth - \varphi)
        \dx \dt\\
&\qquad + \int_{0}^{\tau} \intom \kappa(\vth) \nabla\vth \cdot \nabla \mathcal{T}_{k}(\vth - \varphi) \dx \,dt \leq\int_{0}^{\tau}\intom g \mathcal{T}_{k}(\vth - \varphi) \dx \dt + \intom \mathcal{G}_k (\vth_0 -\varphi(0)) \dx. 
    \end{aligned} 
\end{equation}
\item[{\bf C)}] For every test function $\varphi$ as in {\bf B)}, and for all $k\in \mathbb{N}$ and a.e. $\tau\in (0,T)$, the inequality \eqref{Equi2} holds with equality.
\item[{\bf D)}] There is $\alpha>1/2$ such that the relative energy inequality 
\begin{equation}
\label{ree}
\begin{multlined}
\intom \vth(\tau)-\cQ(\hat\vth)^\alpha\cH_\alpha(\vth(\tau))-\big(\vth_0-\cQ(\hat\vth)^\alpha\cH_\alpha(\vth_0)\big)-\int_0^\tau\intom \bu\cQ(\hat\vth)^\alpha\nabla\cH_\alpha(\vth)\dx\dt \\
+\alpha\int_0^\tau\intom \cQ(\hat\vth)\left|\nabla\left(\frac{\cQ(\vth)}{\cQ(\hat\vth)}\right)\right|^2 \left(\frac{\cQ(\vth)}{\cQ(\hat\vth)}\right)^{-1-\alpha}\dx\dt \leq\int_0^\tau\intom g\left(1-\left(\frac{\cQ(\hat\vth)}{\cQ(\vth)}\right)^\alpha\right)\dx\dt
    \end{multlined} 
\end{equation}
holds for $\tau\in[0,T]$. The functions $\cQ$ and $\cH_\alpha$ are defined as follows\footnote{Recall that $\mu$ is defined in \eqref{muconditions}.}
$$
\forall s>0: \cQ(s)=\int_0^s\kappa(\sigma)\dd \sigma, \quad \cH_\alpha(s)=\int_{\mu}^s\frac1{\cQ(\sigma)^\alpha}\dd\sigma.
$$
\item[{\bf E)}] 
There is $\alpha>1/2$ such that the relative energy equality holds, i.e., \eqref{ree} holds with equality.
\end{itemize}

Finally, any of the statements {\bf A)}--{\bf E)} implies that $\vth\in C([0,T]; L^1(\Omega))$. Consequently, 
{\bf B)}--{\bf C)} are equivalent to the corresponding statements formulated for all $\tau \in (0,T)$.
\end{thm}

This theorem gives several equivalent conditions that imply the time-continuity of the temperature. It shows that $L^1$-continuity follows both from an asymptotic control of the temperature gradient on high level sets, measured by \eqref{Equi1}, and from the validity of the truncated energy inequality \eqref{Equi2} for all admissible test functions. Moreover, this inequality \eqref{Equi2} is equivalent to the corresponding equality. The truncated energy (in)equality \eqref{Equi2}is also equivalent to the relative energy (in)equality expressed in \eqref{ree}. It is also worth noting that the result holds in arbitrary spatial dimensions and requires only low regularity of the velocity field. 

Although the theorem is formulated for strictly positive temperatures, inspired by the Navier--Stokes--Fourier system as in Theorem~\ref{maintheorem}, the statement of Theorem~\ref{maintheoremB} remains valid for general equations not requiring positivity of the solution with a single modification: one cannot consider the entropy equation \eqref{eeB}, but must work solely with the renormalized equation \eqref{iebR}. Indeed, in the course of the proof we do not need the temperature to be nonnegative, except for the step where the entropy is introduced. Therefore, if one considers the problem without the entropy (in)equality, the result remains valid for a general parabolic problem.

\subsection{State of the art}
Foundational works by  Consiglieri \cite{Consiglieri1997, Consiglieri2000} studied both stationary and time-dependent models of shear-thickening fluids with temperature-dependent viscosity and thermal conductivity. While these studies provided rigorous existence results, they did not address the entropy equality, the long-time behaviour, or the continuity of the temperature with respect to time.  More complex existence theory for heat-conducting non-Newtonian fluids was later established in \cite{BuFeMa2009, Bulicek2009}; however, these works did not address additional properties such as time continuity or the validity of energy and/or entropy equality.

From a thermodynamic perspective, Dafermos~\cite{Dafermos1979} introduced the use of entropy as a Lyapunov functional, establishing a link between the Second Law of Thermodynamics and nonlinear stability of fluid systems, see also~\cite{BuMaPr19}.
More recently, Abbatiello et al.~\cite{Abbatiello2022} proved the existence of weak solutions to the three-dimensional Navier--Stokes--Fourier system with non-Newtonian power-law behavior (for $p \geq \frac{11}{5}$) that satisfy the entropy equality. This approach notably avoids the classical DiPerna--Lions renormalization framework \cite{PerenaLions}, by constructing solutions directly fulfilling the entropy equality, even when $p \leq \frac{5}{2}$. Building upon this, the follow-up work \cite{Abbatiello2024} established the nonlinear stability of such entropy-equality solutions by constructing a Lyapunov functional and demonstrating exponential decay to a steady state, even under inhomogeneous boundary temperatures; see also \cite{Dostalik2021} for similar but conditional results.

Considering only the theory for nonlinear heat equations with irregular data (purely integrable or even measure-valued), we overcome the limitations of traditional renormalization approaches. For instance, \cite{Pr97} addressed entropy solutions that are continuous in time, but relied on entropy-like inequalities with multiple cut-off functions, without clearly establishing a renormalized framework or coupling with other equations. The classical existence results for the heat equation with measure data were established by \cite{BoGa89}, while \cite{Po99} introduced renormalized solutions requiring multiple renormalizations in different function spaces. The analysis in \cite{BuCoMa11} extended the existence theory to cases where even the convective term cannot be handled by standard methods, and \cite{AbBuLe24} addressed more complex systems with dissipative heating, albeit without proving continuity with respect to time.

The present article contributes to this line of research by focusing on a two-dimensional non-Newtonian heat-conducting fluid and establishing the existence of entropy-equality solutions. The main novelty of our work is the time regularity of the solution: we show that the temperature is continuous in time, which allows us to avoid classical renormalization procedures altogether. Remarkably, in most works on heat-conducting fluids, such issues are usually omitted, and even the attainment of the initial condition in the strong topology is typically not guaranteed. In this paper, all the above-mentioned drawbacks are resolved, and the validity of the entropy equation appears to be the key step allowing us to avoid the classical renormalized solution notions while still proving the desired strong results.

\section{Proof of the existence of an entropy solution - Theorem~\ref{maintheorem}}\label{sec:3}
In the first four sections, following the methods of \cite{BuFeMa2009}, we prove the existence of a weak solution satisfying \eqref{me} and \eqref{ieb}. The proof requires modification due to the different character of the boundary conditions. We then show that the initial condition for the temperature is attained in the sense of \eqref{sqrtconvergenceth0A}.  

Next, we prove that the solution satisfies \eqref{ee}, adopting the methods from \cite{Abbatiello2022}. Finally, we show that the temperature $\vth$ is continuous in $L^1$ with respect to time, and hence the initial condition is attained in the sense of Theorem~\ref{maintheorem}.

\subsection{Galerkin approximations}
We now define the Galerkin approximation. First, we denote
$$
W^{2,2}_{0,\diver}\om := \overline{\{\bu \in (C_0^{\infty}(\Omega))^2; \diver \bu = 0\}}^{\|\cdot\|_{2,2}}.
$$
Then, following \cite[Appendix, Theorem 4.11]{MaNeRoRu95}, we obtain the existence of a basis $\{\bw_j\}_{j=1}^\infty$ that is orthogonal in $W^{2,2}_{0,\diver}\om$ and orthonormal in $L^2(\Omega)$. Note that $\norm{\bw_j}{1,q}{ }$ is bounded for all $j \in \nn$, $q \in [1,\infty)$.  
Similarly, we can find a sequence $\{w_j\}_{j=1}^\infty$ forming an orthonormal basis of $L^{2}\om$, which is orthogonal in $W^{1,2}_{0}\om$.

Then, for fixed $N,M\in\nn$, we find 
\begin{align*}
    \bu^{N,M}(x,t) &= \summing{i}{N} c_i^{N,M}(t) \bw_i(x), \\
    \vth^{N,M}(x,t) &= \left(\summing{i}{M} d_i^{N,M}(t) w_i(x)\right) + \htet(x), 
\end{align*}
such that $(\pmb{c}^{N, M}, \pmb{d}^{N, M}):=(c_1^{N,M},\dots,c_N^{N,M}, d_1^{N,M},\dots,d_M^{N,M}):(0,T^\ast)\to\rr^{N+M}$, for $T^\ast \in (0,T)$, is the maximal solution to the following system of ordinary differential equations
\begin{equation}\label{ODRgal1}
    \begin{aligned}
       &\intom  \dert \bu^{N,M}\bw_j\dx + \int_\Omega \tens^\ast\left(\vth^{N,M}, D\bu^{N,M}\right):D\bw_j \dx = \\
    & \quad \quad \quad \quad  \quad \quad \quad \quad = \intom (\bu^{N,M} \otimes \bu^{N,M}):\nabla(\bw_j) \dx  +  \dual{\beh}{\bw_j}, \\
   &\intom \dert \vth^{N,M} w_k \dx + \intom \kappa(\vth^{N,M}) \nabla\vth^{N,M} \cdot (\nabla w_k)  \dx = \\
    & = \intom \tens^\ast\left(\vth^{N,M}, D\bu^{N,M}\right):D\bu^{N,M} w_k \dx + \int_\Omega \vth^{N,M} \bu^{N,M} \cdot (\nabla w_k) \dx, 
    \end{aligned} 
\end{equation}
where $j\in \{1, \dots, N\}$ and $ k\in \{1, \dots, M\}$.

We equip the above system of ordinary differential equations \eqref{ODRgal1} with the following initial conditions:
\begin{equation}\label{ODRic}
    \begin{aligned}
        \bu^{N,M}(x,0) &= \bu_0^{N}(x)  &&\textrm{on $\Omega$},
        \\
        \vth^{N,M}(x,0) &= \vth_0^{N,M}(x)  &&\textrm{on $\Omega$},
    \end{aligned}
\end{equation}
where $\bu_0^{N}(x) = \sum_{i = 1}^N c_{i,0}^{N} \bw_i(x)$ is the projection of $\bu_0$ onto the linear hull of $\{\bw_j\}_{j=0}^N$ in $L^2(\Omega)$. The approximative initial conditions for the temperature  $\vth_0^{N,M}$ is defined in the following manner. We define $\vth_0^N(x):=\min\{\vth_0(x),N\}$ 
and then  define $\vth_0^{N,M}$ as a projection of $\vth_0^{N}$ onto the linear hull of $\{w_j\}_{j=0}^M$ in $L^2(\Omega)$, hence $\vth_0^{N,M} = \sum_{j=1}^M d_{0,j}^{N,M} w_j$. Note, that
\begin{equation}\label{icconvergences}
    \begin{aligned}
          \vth_0^{N,M} &\xrightarrow {\text{$M \to \infty$}} \vth_0^{N} &&\textrm{in $L^2(\Omega)$},\\ 
          \bu_0^{N} &\xrightarrow{\text{$N \to \infty$}} \bu_0 &&\textrm{in $L^2(\Omega;\rr^2)$,}\\ 
          \vth_0^{N} &\xrightarrow{\text{$N \to \infty$}} \vth_0 &&\textrm{in $L^1(\Omega)$}, 
    \end{aligned}
\end{equation}
where the first two convergence results hold due to the completeness of an orthonormal basis in Hilbert spaces.

By Carathéodory's existence theory, it is not difficult to obtain a maximal solution to the system \eqref{ODRgal1}--\eqref{ODRic} on a time interval $(0,T^*)$. Showing that $\bu^{N,M}$ and $\vth^{N,M}$ are bounded with respect to the time variable then implies that we can set $T^* = T$.

\subsection{Limit passage  \texorpdfstring{$M\to \infty$}{M}} 
To pass to the limit as $M \to \infty$ and to establish existence on the time interval $(0,T)$, we need to derive estimates for $\bu^{N,M}$ and $\vth^{N,M}$ that are uniform with respect to $M$. For the velocity estimates, we follow the classical procedure. We multiply the $j$-th equation in \eqref{ODRgal1}$_1$ by $c_j^{N,M}$, sum the results over $j=1,\ldots, N$, and use the fact that $\diver \bu^{N,M}=0$ to obtain
\begin{equation}
\frac12 \frac{\rm{d}}{\dt} \|\bu^{N,M}(t)\|_2^2 + \int_{\Omega} \tens (\vth^{N,M},D\bu^{N,M}) : D\bu^{N,M}\dx = \dual{\beh}{\bu^{N,M}}. \label{NM1}
\end{equation}
Then, using the assumption \eqref{tensor}, the assumptions on the data \eqref{conditions}, the Korn inequality \eqref{korn}, the properties of the approximation of the initial data \eqref{icconvergences}, the interpolation \eqref{galnir}, and the identity \eqref{ODRgal1}$_1$ for estimating the time derivative, we deduce the uniform bounds
\begin{equation}\label{estimateu}
    \begin{aligned}
        \norm{\bu^{N,M}}{L^\infty(L^2)}{} + \norm{\bu^{N,M}}{L^p(V)}{} + \norm{\bu^{N,M}}{L^{2p}(Q)}{} + \norm{\dert \bu^{N,M}}{L^{p'}((W^{2,2}_{0,\diver})^\ast)}{} \leq C.
    \end{aligned} 
\end{equation}
These estimates are crucial for passing to the limit and establishing existence on the interval $(0,T)$.

Next, we focus on the temperature estimates. We multiply the $(M+k)$-th equation in \eqref{ODRgal1} by $d_k^{N,M}$ and sum over $k = 1, \dots , M$. This yields (using the fact that $\diver  \bu^{N,M}=0$) the identity
$$
\begin{aligned}
&\frac12 \frac{\rm{d}}{\dt} \|\vth^{N,M}(t)-\htet\|_2^2 + \intom \kappa(\vth^{N,M}) \nabla\vth^{N,M} \cdot \nabla(\vth^{N,M}-\htet)  \dx = \\
    &\qquad  = \intom \tens^\ast\left(\vth^{N,M}, D\bu^{N,M}\right):D\bu^{N,M}(\vth^{N,M}-\htet) \dx + \int_\Omega  \htet\, \bu^{N,M} \cdot \nabla\vth^{N,M}\dx.
\end{aligned}
$$
Then, using \eqref{kap} and \eqref{tensor}, the Young inequality, and the Poincar\'e inequality, while integrating over the time interval $(0,t)$ for any $t \in (0,T)$, we have
\begin{equation}\label{estimateth0}
    \begin{aligned}
         &\|(\vth^{N,M} -  \hat{\vth})(t)\|_2^2 +  \int_0^t \norm{\vth^{N,M} - \hat{\vth}}{1,2}{2} \, \mathrm{d}\tau \\
         & \qquad \leq C\|(\vth^{N,M} -  \hat{\vth})(0)\|_2^2  + C \int_0^t \Big( 1 + \norm{D\bu^{N,M}}{2p}{2p} + \norm{\hat{\vth}}{\infty}{2} \norm{\bu^{N,M}}{2}{2} + \norm{\hat{\vth}}{2}{2} \Big) \, \mathrm{d}\tau \\
         &\qquad \leq C \left(1+ \int_0^t \norm{D\bu^{N,M}}{2p}{2p} \, \mathrm{d}\tau\right).
    \end{aligned} 
\end{equation}
By the choice of the basis, we have
\begin{align}\label{Dup}
    \int_0^t\norm{D\bu^{N,M}}{2p}{2p}\dtau \leq \sup_{j \leq N} \norm{D\bw_j}{2p}{2p} \int_0^t\left(\sum_{j=1}^N{c}_j^{N,M}\right)^{2p}\dtau  \overset{\eqref{estimateu}}\leq C(N).
\end{align}
Consequently, it follows from \eqref{estimateth0} that 
\begin{equation}\label{estimateth}
    \begin{aligned}
   \sup_{t\in(0,T)}\|(\vth^{N,M} -  \hat{\vth})(t)\|_2^2 + \int_0^{T} \norm{\vth^{N,M} - \hat{\vth}}{1,2}{2} \dt  \leq C(N).
    \end{aligned} 
\end{equation}
From \eqref{ODRgal1}, we deduce
\begin{equation}\label{estimatedth}
    \begin{aligned}
   &\norm{\dert \vth^{N,M}}{(L^2(W_0^{1,2})^\ast)}{ } \leq  \intt C(1+\norm{D\bu^{N,M}}{2p}{2p} + \norm{\vth^{N,M}}{4}{4}+\norm{\bu^{N,M}}{4}{4} +  \normg{\vth^{N,M}}{2}{2})\dt \leq 
    C(N).
    \end{aligned} 
\end{equation}
For any fixed $N \in \nn$, due to the $M$-independent estimates \eqref{estimateu}--\eqref{estimatedth}, we can find non-relabeled subsequences $\{\bu^{N,M}\}_{M=1}^\infty, \{\vth^{N,M}\}_{M=1}^\infty$ such that, as $M \to \infty$, we obtain the following convergences
\begin{align}
    \dert \bu^{N,M} &\rightharpoonup \dert \bu^{N} && \text{in $L^{p'}(0,T;(W^{2,2}_{0,\diver})^\ast)$},\label{dertu}\\
    \bu^{N,M} &\rightharpoonup^\ast \bu^{N} && \text{in $L^{\infty}(0,T;L^2_{0,\diver})$},\\
    \bu^{N,M} &\rightharpoonup \bu^{N} && \text{in $L^{p}(0,T;V)$},\\
    \dert \vth^{N,M} &\rightharpoonup \dert \vth^{N} && \text{in $L^{2}(0,T;(W^{1,2}_0)^\ast)$},\label{timeMB}\\ 
    \label{convergenceofdtheta}
    \vth^{N,M} &\rightharpoonup^\ast \vth^{N} && \text{in $L^{\infty}(0,T;L^2\om)$},\\
    \vth^{N,M} &\rightharpoonup \vth^{N} && \text{in $L^{2}(0,T;W^{1,2}\om)$}. \label{convergenceoftheta}
\end{align}
With the use of the Aubin--Lions lemma (see \cite[Corollary~9]{Si1987} and also \cite{Ro90}), Lebesgue-Sobolev interpolation \eqref{galnir}, and the fact that all norms are equivalent on finite dimensional spaces, we can additionally obtain strong convergence results
\begin{align}
         & \vth^{N,M} \to \vth^{N} && \text{in $L^{r}(Q)$ for $r \in \left[1,4\right)$},\label{thstrong}\\
    & \bu^{N,M} \to \bu^{N} && \text{in $L^{r}(0,T;W^{1,r}\om)$ for $r \in \left[1,+\infty\right)$}.\label{ustrong}
\end{align}
We can now pass to the limit $M\to\infty$ in \eqref{ODRgal1} and obtain
\begin{equation}\label{1gal}
    \begin{aligned}
       \dual{\dert \bu^{N}}{\bw_j}_{W^{2,2}_{0,\diver}} + &\int_\Omega \tens^\ast\left(\vth^{N}, D\bu^{N}\right):D\bw_j \,dx = &&\\
    & \quad \quad \quad = \intom (\bu^{N} \otimes \bu^{N}):\nabla\bw_j \,dx  +  \dual{\beh}{\bw_j}, 
 && \text{for all $j \in \{1, 2, \dots, N\}$}
     \end{aligned} 
\end{equation}
\begin{equation}\label{2gal}
    \begin{aligned}
   \dual{\dert \vth^{N}}{\psi}_{W^{1,2}_0} + &\intom \kappa(\vth^{N}) \nabla\vth^{N} \cdot \nabla \psi   \dx = &&\\
    &= \intom \tens^\ast\left(\vth^{N}, D\bu^{N}\right):D\bu^{N} \psi  \dx + \int_\Omega \vth^{N} \bu^{N} \cdot \nabla \psi  \dx, && \forall \psi \in W_0^{1,2}(\Omega).
    \end{aligned} 
\end{equation}
Furthermore, using the classical parabolic embedding, the fact that $\htet$ is time independent, and the convergence results \eqref{timeMB} and \eqref{convergenceoftheta}, we have
\begin{align}\label{thetaNcontinuity}
     \vth^N \in C([0,T];L^{2}).
 \end{align}

\subsection{Minimum principle}
We show that $\vth^N \geq \mu$ for almost all $(t,x) \in (0,T)\times \Omega$, where $\mu$ is defined in \eqref{muconditions}. To this end, we set
\begin{align*}
    \psi(\tau,x) = \min\{0, \vth^N(\tau, x)-\mu\} = (\vth^N(\tau, x)-\mu)_{-} \leq 0
\end{align*}
in \eqref{2gal} and integrate the resulting expression until $t \in (0,T)$. Note that $\psi$ has a zero trace due to the choice of $\mu$. 
Using $\diver \bu^{N}=0$, we have
\begin{equation}\label{minimumest}
    \begin{aligned}
   &\ \int_0^t \frac12\dert\norm{\psi}{2}{2}\dtau + \int_0^t \intom \kappa(\vth^{N}) |\nabla \psi|^2 \dx \dtau  = \int_0^t \intom  \tens^\ast\left(\vth^{N}, D\bu^{N}\right):D\bu^{N} \psi  \dx\dtau \le 0,
    \end{aligned} 
\end{equation}
where we used the nonnegativity of\footnote{It follows from \eqref{tensor} by using the following inequality
$$
0 \le (\tens^\ast(\vth^{N}, D\bu^{N}) - \tens^\ast(\vth^{N}, \bO)) : (D\bu^{N}-\bO) = \tens^\ast(\vth^{N}, D\bu^{N}):D\bu^{N}.
$$
}
$\tens^\ast (\vth^{N}, D\bu^{N}):D\bu^{N}$ and the nonpositivity of $\psi$.
This implies 
\begin{align*}
& 0 \geq \int_0^t \dert \intom \psi^2 \dx \dtau = \norm{(\vth^N(t)-\mu)_{-}}{2}{2} - \norm{(\vth^N_0-\mu)_{-}}{2}{2}=\norm{(\vth^N(t)-\mu)_{-}}{2}{2},
\end{align*}
which yields 
\begin{align}\label{minp}
    &\vth^N \geq \mu &&\text{for a. a. } (t,x) \in (0,T)\times \Omega.
\end{align}

\bigskip

To summarize, we found functions 
\begin{align*}
    &\bu^N:Q \to \rr^2 &\vth^N:Q\to\rr
\end{align*}
such that 
\begin{equation}\label{propertiesN}
    \begin{aligned}
    &\bu^N\in L^\infty(0,T;L^2_{0,\diver})\cap L^p(0,T;V), && \dert\bu^N\in L^{p'}(0,T;(W^{2,2}_{0,\diver})^\ast),\\
    &\vth^N\in C([0,T];L^2\om)\cap L^2(0,T;W^{1,2}\om), && \dert\vth^N\in L^{2}(0,T;(W^{1,2}_0)^*),\\
    & \vth^N \geq \mu \text{ a.a. on }Q.
\end{aligned}
\end{equation}
Furthermore, $\bu^N$ and $\vth^N$ fulfill the system \eqref{1gal}--\eqref{2gal}, and satisfy the initial conditions
\begin{align*}
    &\bu^{N}(0) = \bu^{N}_0 \quad \textrm{and} \quad \vth^{N}(0) = \vth^{N}_0.
\end{align*} 
We omit the discussion of the attainment of the initial conditions here, but we focus on it in the next section, where a more delicate procedure must be employed.

\subsection{Limit passage  \texorpdfstring{$N\to \infty$}{N}}
\label{sec:ntn}
We provide estimates on $\bu^N$ and $\vth^N$ that are uniform with respect to $N$ and obtain weakly convergent subsequences.
Since the estimates for $\bu^{N,M}$ were uniform with respect to $N$, we can use weak lower semicontinuity and obtain from \eqref{estimateu} that $\bu^N$ also satisfies
\begin{equation}\label{estimateu2}
    \begin{aligned}
        \sup_{t \in (0,T)} \|\bu^{N}(t)\|_2 + \norm{\bu^{N}}{L^{p}(V)}{}+ \norm{\bu^{N}}{L^{2p}(Q)}{} +  \norm{\dert \bu^N}{L^{p'}((W^{2,2}_{0,\diver})^\ast)}{} \leq C.
    \end{aligned} 
\end{equation}
In addition, similarly to \eqref{NM1}, we obtain the following identity 
\begin{equation}
\frac12 \frac{\rm{d}}{\dt} \|\bu^{N}(t)\|_2^2 + \int_{\Omega} \tens (\vth^{N},D\bu^{N}) : D\bu^{N}\dx = \dual{\beh}{\bu^{N}}. \label{NM1N}
\end{equation}

\bigskip

The same procedure for $\vth^N$ cannot be used since \eqref{estimateth} is not uniform with respect to $N$. Therefore, we proceed differently. 
We consider arbitrary $\alpha \in (0,1)$, denote $K:= \|\hat{\vth}\|_{\infty}$ and set
\begin{align*}
    &\psi := 1 - \left(\frac{K}{K + \vth^N - \hat{\vth}}\right)^\alpha
\end{align*}
as a test function in \eqref{2gal}. Note that $\psi \in W^{1,2}_0(\Omega) \cap L^{\infty}(\Omega)$ due to the properties of $\hat{\vth}$ and the minimum principle \eqref{minp}. Integrating the result over $t \in (0,T)$, using the positivity of $\tens^N:D\bu^{N}$ and the fact $\psi \leq 1$, we have
\begin{equation}\label{testN}
    \begin{aligned}
   & \int_0^t \dual{\dert \vth^{N}}{1 - \left(\frac{K}{K + \vth^N-\hat{\vth}}\right)^\alpha}_{W^{1,2}_0}\dtau  -
  \int_0^t \intom \left(\kappa(\vth^{N}) \nabla\vth^{N} - \vth^{N} \bu^{N}\right) \cdot \nabla \left(\frac{K}{K + \vth^N-\hat{\vth}}\right)^\alpha \dtau  \leq \\
    & \leq \int_0^t \intom \tens^N:D\bu^{N} \dx \dtau \leq C.
    \end{aligned} 
\end{equation}
The last inequality holds by \eqref{estimateu2}, and the third part of \eqref{tensor}.
We define $H^\alpha: (0,\infty)^2 \to \rr$ as
\begin{align*}
    H^\alpha(s,\sigma) = \int_0^s \left(\frac{K}{K + \tau-\sigma}\right)^\alpha\dtau
\end{align*}
and we observe that there exist constants $c_1, c_2 > 0$ such that for all $s \geq \mu$ and $0 \le \sigma \le K$, it holds
\begin{align}\label{halpha}
   c_1(\alpha)s^{1-\alpha} \leq H^\alpha(s,\sigma) \leq c_2(\alpha)s^{1-\alpha}. 
\end{align}
Next, we evaluate all terms in \eqref{testN}. For the first term, we have  
\begin{align}
        \dual{\dert \vth^{N}}{1 - \left(\frac{K}{K + \vth^N-\hat{\vth}}\right)^\alpha}_{W^{1,2}_0} = \frac{\rm{d}}{\dt}\, \norm{\vth^{N} - H^\alpha(\vth^N, \htet))}{1}{ }.
\end{align}
We now want to simplify the term in \eqref{testN}, which contains heat conductivity $\kappa$. Using the Young inequality, we have 
\begin{equation}\label{who}
    \begin{aligned}
     -\kappa(\vth^{N}) \nabla\vth^{N} \cdot \nabla \left(\frac{K}{K + \vth^N-\hat{\vth}}\right)^\alpha &= \frac{\alpha\kappa(\vth^{N})K^\alpha}{\left(K + \vth^N-\hat{\vth}\right)^{\alpha+ 1}} |\nabla \vth^N|^2   - \frac{\alpha\kappa(\vth^{N})K^\alpha}{\left(K + \vth^N-\hat{\vth}\right)^{\alpha+ 1}} \nabla \vth^N\cdot  \nabla\hat{\vth}\\
     &\ge \frac{\alpha\kappa(\vth^{N})K^\alpha}{2\left(K + \vth^N-\hat{\vth}\right)^{\alpha+ 1}} |\nabla \vth^N|^2   - \frac{\alpha\kappa(\vth^{N})K^\alpha}{2\left(K + \vth^N-\hat{\vth}\right)^{\alpha+ 1}} |\nabla\hat{\vth}|^2\\
     &\ge \frac{\alpha\kappa(\vth^{N})K^\alpha}{2\left(K + \vth^N-\hat{\vth}\right)^{\alpha+ 1}} |\nabla \vth^N|^2   - C(\alpha,K,\mu,\overline{\kappa})|\nabla\hat{\vth}|^2,
    \end{aligned}
\end{equation}
where we used \eqref{kap}, the upper bound on $\kappa$. 
 Using $\diver \bu^N = 0$, we estimate the last term on the left-hand side of \eqref{testN} as 
\begin{align*}
   -\int_0^t \intom  \vth^{N} \bu^{N} \cdot \nabla \left(\frac{K}{K + \vth^N-\hat{\vth}}\right)^\alpha \dx \dtau 
   &= -\int_0^t \intom  \hat{\vth} \bu^{N} \cdot \nabla \left(\frac{K}{K + \vth^N-\hat{\vth}}\right)^\alpha \dx \dtau   \\
   &= \int_0^t \intom  \nabla \hat{\vth} \cdot \bu^{N}  \left(\frac{K}{K + \vth^N-\hat{\vth}}\right)^\alpha \dx \dtau \\
   &\le C(K,\alpha,\mu) \int_0^t \intom  |\nabla \hat{\vth}|^2 +  |\bu^{N}|^2 \dx \dtau \le C(K,\alpha,\mu),
\end{align*}
where the last inequality follows from the assumption \eqref{conditionstheta} on $\hat{\vth}$ and the uniform estimate \eqref{estimateu2}.
%
%
%
Putting all the estimates together, we can rewrite \eqref{testN} as 
\begin{equation*}
    \begin{aligned}
        \norm{\vth^{N}(t) - H^\alpha(\vth^N(t), \htet))}{1}{ } + \int_0^t \intom\frac{\alpha\kappa(\vth^{N})K^\alpha |\nabla \vth^N|^2}{2\left(K + \vth^N-\hat{\vth}\right)^{\alpha+ 1}}  \dtom \leq \norm{\vth^{N}_0 - H^\alpha(\vth^N_0, \htet))}{1}{ } +C(\alpha,\mu,K,\overline{\kappa}).
    \end{aligned}
\end{equation*}
 Using the minimum principle \eqref{minp}, the assumption \eqref{kap}, the assumptions on the initial data \eqref{conditions} and \eqref{icconvergences}, and \eqref{halpha} bounds on $H^\alpha$, we deduce from the above estimate that 

%
%
%
\begin{equation}\label{estimatecombofinal}
    \begin{aligned}
        &\sup_{t \in (0,T)} \norm{\vth^{N} - H^\alpha(\vth^N, \htet))}{1}{ } + \intq \frac{|\nabla \vth^N|^2}{(\vth^N)^{\alpha +1}} \dq \leq C.
    \end{aligned}
\end{equation}
Next, using \eqref{halpha} and the Young inequality, we obtain the estimate
\begin{align*}
    \intom |\vth^{N}| \dx 
    &\leq \norm{\vth^{N} - H^\alpha(\vth^N, \htet)}{1}{ } + \intom c_2 (\vth^{N})^{1-\alpha} \dx\le \norm{\vth^{N} - H^\alpha(\vth^N, \htet)}{1}{ } + C + \frac{1}{2}\intom |\vth^{N}| \dx.
\end{align*}
Consequently, \eqref{estimatecombofinal} implies
\begin{align}\label{estimatethn}
    \norm{\vth^N}{L^\infty(L^1)}{ } \leq C.
\end{align}
%
%
We can now use \eqref{estimatecombofinal} together with \eqref{conditionstheta}--\eqref{muconditions} to estimate
\begin{align*}
    &\intq \left|\nabla \left((\vth^{N})^{\frac{1-\alpha}{2}}-\hat{\vth}^{\frac{1-\alpha}{2}}\right)\right|^2 \dq 
    \leq C\intq \frac{\left|\nabla \vth^{N}\right|^2}{\left(\vth^N\right)^{\alpha+1}}  
    + \frac{|\nabla \hat{\vth}|^2}{(\hat{\vth})^{\alpha+1}} \, \dq 
    \le C,
\end{align*}
and using the Poincar\'{e} inequality and the classical Sobolev embedding, we conclude that for all finite $s \in [1,\infty)$ it holds
\begin{align*}
     \norm{(\vth^{N})^{\frac{1-\alpha}{2}}-(\hat{\vth})^{\frac{1-\alpha}{2}}}{L^2(L^s)}{}\leq C.
\end{align*}
Furthermore, it follows from the boundedness of $\hat{\vth}$ that for all $s<\infty$
\begin{align}\label{estimatethn2}
    \norm{(\vth^{N})^{\frac{1-\alpha}{2}}}{ L^2(L^s)}{}\leq C.
\end{align}
Interpolating \eqref{estimatethn} and \eqref{estimatethn2}, and allowing $s$ to be arbitrarily large and $\alpha>0$ arbitrarily small, we have the estimate
\begin{align}\label{estimatethN}
    &\norm{\vth^N}{L^r(Q)}{ } \leq C &&\text{for } r \in [1, 2).
\end{align}

Let us now estimate the gradient of $\vth^N$ using Hölder's inequality:
\begin{align*}
    \intq | \nabla \vth^N |^t \, \dq 
    &= \intq \left(\frac{| \nabla \vth^N |}{\left(\vth^N\right)^\frac{\alpha+1}{2}} \right)^t \left({\vth^N}\right)^\frac{t(\alpha+1)}{2} \, \dq \\
    &\leq \left(\intq \frac{| \nabla \vth^N |^2}{\left(\vth^N\right)^{\alpha+1}} \, \dq \right)^\frac{t}{2} 
    \left(\intq \left({\vth^N}\right)^\frac{t(\alpha+1)}{2-t} \, \dq \right)^\frac{2-t}{2},
\end{align*}
where the first factor is bounded by \eqref{estimatecombofinal} and the second factor is finite by \eqref{estimatethN} if $t < \frac{4}{3 + \alpha}$.
Hence, taking $\alpha >0$ arbitrarily small, we obtain
\begin{align}\label{estimategradthN}
    \norm{\nabla \vth^N}{L^t(Q)}{} \leq C && \text{for } t \in \left[1, \frac{4}{3}\right).
\end{align}

Finally, we focus on the estimate for the time derivative of $\vth^N$. Take $K:=\{\psi \in W^{1,5}_0(\Omega); \; \|\psi\|_{1,5} \le 1\}$. Using embedding $W^{1,5}\hookrightarrow L^{\infty}$, H\"{o}lder's inequality, Young inequality, and assumption \eqref{kap}, it follows from \eqref{2gal} that for almost all time $t\in (0,T)$ we have 
$$
\begin{aligned}
\|\partial_t \vth^N (t)\|_{(W^{1,5}_0)^*} &= \sup_{\psi \in K} \dual{\partial_t \vth^{N}(t)}{\psi} \\
&= \sup_{\psi \in K} \left(-\intom \kappa(\vth^{N}) \nabla\vth^{N} \cdot \nabla \psi   \dx + \intom \tens^N(t):D\bu^{N} \psi  \dx + \int_\Omega \vth^{N} \bu^{N} \cdot \nabla \psi  \dx \right)\\
&\le \overline{\kappa}\|\nabla\vth^{N}(t)\|_{\frac{5}{4}} +C \intom \tens^N(t):D\bu^{N}(t)  \dx + \|\vth^{N}(t) \, \bu^{N}(t)\|_{\frac{5}{4}}\\
&\le C \intom 1 + |\nabla\vth^{N}(t)|^{\frac{5}{4}} + \tens^N(t):D\bu^{N}(t)  + |\vth^{N}(t)|^{\frac{20}{11}} +  |\bu^{N}(t)|^4 \dx .
\end{aligned}
$$
Integrating the above inequality over the time interval $(0,T)$, and using the uniform estimates 
\eqref{estimatethN}, \eqref{estimategradthN}, and \eqref{estimateu2}, we obtain
%
\begin{align}\label{estimatedthN}
    &\norm{\dert \vth^N}{L^1\left(\left(W^{1,5}_0\right)^\ast\right)}{ }\leq C.
\end{align}

The uniform estimates \eqref{estimateu2}, \eqref{estimatecombofinal}, \eqref{estimatethn}, \eqref{estimatethN}, and \eqref{estimategradthN},  together with the Aubin--Lions lemma, provide us with the following convergence results for a subsequence that we do not relabel:
\begin{align}
    \dert \bu^{N} &\rightharpoonup \dert \bu && \text{in $L^{p'}(0,T;(W^{2,2}_{0,\diver})^\ast)$}, \label{Nconvergence1}\\
    \bu^{N} &\rightharpoonup^\ast \bu && \text{in $L^{\infty}(0,T;L^2_{0,\diver})$},\label{Nconvergence2}\\
    \bu^{N} &\rightharpoonup \bu && \text{in $L^{p}(0,T;V)$},\label{Nconvergence3}\\
     \bu^{N} &\to \bu && \text{in $L^{r}(0,T;L^r)$ for $r \in \left[1,{2}p\right)$}, \label{Nconvergence4}\\
     \tens^{N} &\rightharpoonup \tens && \text{in $L^{p'}(0,T;L^{p'}(\Omega;\rr^{2\times2}_{sym}))$}.\label{Nconvergence5}\\
 \vth^{N} &\rightharpoonup \vth && \text{in $L^{s}(0,T;L^s)$ for $s \in \left[1,2\right)$},\label{Nconvergence6}\\
     \nabla \vth^{N} &\rightharpoonup \nabla \vth && \text{in $L^{t}(0,T;L^t)$ for $t \in \left[1,{4}/{3}\right)$},\label{Nconvergence7}\\
    (\vth^{N})^\alpha &\rightharpoonup  (\vth )^\alpha && \text{in $L^{2}(0,T;W^{1,2}\om)$ for $\alpha \in \left[0,{1}/{2}\right)$}. \label{Nconvergence8}
\end{align}
Furthermore, the generalized version of the Aubin--Lions lemma (see \cite{Ro90}), together with the uniform bound~\eqref{estimatedthN} and the above convergence results, implies
\begin{align}
     &\vth^{N} \to \vth && \text{almost everywhere in $Q$},\label{aeNconvergence}\\
\label{strongthetaNconvergence}
     &\vth^{N} \to \vth && \text{in $L^{s}(0,T,L^{s})$ for $s\in \left[1,2\right)$}.
\end{align}
Consequently, from \eqref{estimatethn} and \eqref{strongthetaNconvergence}, we can also deduce that
\begin{align}\label{LinftyL1}
    \vth \in L^\infty(0,T;L^1\om).
\end{align}

We now test equation \eqref{1gal} by $\varphi \in \mathcal{D}(0,T)$ and take the limit $N \to \infty$. Convergences \eqref{Nconvergence1}-\eqref{Nconvergence5} imply 
\begin{equation}\label{product}
    \begin{aligned}
       &\intt\dual{\dert \bu}{\pmb{\psi}}_{W_{0,\diver}^{2,2}}\varphi\dt + \intq \varphi \tens:D\pmb{\psi}\dq =\intq\varphi(\bu \otimes \bu):\nabla\pmb{\psi}  \dq   +   \intt\dual{\beh}{\pmb{\psi} }_V\varphi\dt, 
     \end{aligned} 
\end{equation}
for all  $\pmb{\psi} \in W_{0,\diver}^{2,2}(\Omega)$.
We will now show that $\dert\bu \in L^{p'}(0,T;V^\ast)$. We know that by the definition of the weak time derivative 
\begin{align}\label{weakdu}
    & \intt\dual{\dert \bu}{\pmb{\psi}}_{W_{0,\diver}^{2,2}}\varphi\dt= - \intq \bu \dert(\pmb{\psi} \varphi) \dq =  \mathcal{F}(\pmb{\psi} \varphi),&& \forall \varphi \in \mathcal{D}(0,T), \forall \pmb{\psi} \in W^{2,2}_{0,\diver}(\Omega)
\end{align}
where 
\begin{align*}
    \mathcal{F}(\pmb{\psi}\varphi) := \intq \varphi (\bu \otimes \bu):\nabla\pmb{\psi} - \varphi\tens:D\pmb{\psi} \dq   +   \intt \dual{\beh}{\varphi \pmb{\psi} }_V\dt.
\end{align*}
Furthermore, $\mathcal{F}$ is bounded for all $\bw \in L^p(0,T;V)$, since 
\begin{align*}
    &|\mathcal{F}(\bw)|\leq \norm{\bu}{L^{2p'}(Q)}{2}\norm{\nabla \bw}{L^p(Q)}{ }+ \norm{\tens}{L^{p'}(Q)}{}\norm{D\bw}{L^p(Q)}{ }+ \norm{\beh}{L^{p'}(V^\ast)}{}\norm{\bw}{L^p(V)}{ }.
\end{align*}
We hence have $\mathcal{F} \in L^{p'}(0,T;V^\ast)$. Since the set $span\{\varphi\pmb{\psi} ;\varphi \in \mathcal{D}(0,T), \pmb{\psi} \in W^{2,2}_{0,\diver}(\Omega)\}$ is dense in $L^p(0,T;V)$, there exists a uniquely defined extension of $\dert \bu$ (denoted again $\dert \bu$) such that 
\begin{align} \label{dertub}
    \dert \bu \in L^{p'}(0,T;V^\ast).
\end{align}
This extension is defined by $\dual{\dert \bu}{\bw}_{L^p(V)} =  \mathcal{F}(\bw)$ for all $\bw \in L^{p}(0,T;V)$.
Hence, we can conclude that for every $\bw\in L^p(0,T,V)$ it holds
\begin{equation}\label{limitequation}
  \begin{aligned}
       &\intt\dual{\dert \bu}{ \bw}_V\dt + \intq \tens:D\bw \dq =\intq (\bu \otimes \bu):\nabla\bw \dq   +   \intt\dual{\beh}{\bw}_V\dt. 
\end{aligned}   
\end{equation}
Additionally, by \eqref{Nconvergence3} and \eqref{dertub}, we have $\bu \in C([0,T];L_{0,\diver}^2)$, and the initial condition $\bu_0$ is attained in $L^2(\Omega)$. 

It remains to identify $\tens$ with $\tens^\ast\left(\vth, D\bu\right)$ by the Minty trick. 
Note that if we extended $f$  by zero also for times $t>T$, the construction of the solution could have been done on the interval $(0,T^*)$ with arbitrary $T^*>T$. The identity \eqref{limitequation} then remains valid with $T$ replaced by $T^*$, and all quantities are well defined in $(0,T^*)$, in particular $\bu \in C([0,T^*];L_{0,\diver}^2)$. Setting $\bw:=\bu \chi_{[0,T]}$ and using $\diver \bu=0$, we obtain
\begin{equation}\label{limitequationE}
  \begin{aligned}
       &\|\bu(T)\|_2^2 + 2\int_0^{T} \intom \tens:D\bu \dq =  \|\bu_0\|_2^2 +  2\int_0^{T} \dual{\beh}{\bu}_V\dt. 
\end{aligned}   
\end{equation}
We recall that $\tens^N := \tens^\ast(\vth^N, D\bu^{N})$. Let us take any $\dens \in L^p(Q;\rr^{2\times2}_{sym})$ and extend it by zero for $t>T$. It holds for $\tau\in(0,T^*-T)$
\begin{equation}\label{MintyN}
    \begin{aligned}
    0\stackrel{\eqref{tensor}_1}{\leq} &\intq (\tens^{N} - \tens^\ast(\vth^N, \dens )):(D\bu^{N} - \dens) \dq \le \int_0^{T+\tau}\intom (\tens^{N} - \tens^\ast(\vth^N, \dens )):(D\bu^{N} - \dens) \dq \\
    &= \frac{1}{2}\left(\|\bu^N_0\|_2^2 - \|\bu^N(T+\tau)\|_2^2\right) \\
    &+ \int_0^{T} \dual{\beh}{\bu^{N}}_V \dt - \int_0^{T}\intom \tens^{N}: \dens   +  \tens^\ast(\vth^N, \dens ):(D\bu^{N} - \dens ) \dq,
\end{aligned} 
\end{equation}
where for the last equality we used the identity \eqref{NM1N} integrated over $(0,T+\tau)$.
Next, we integrate the above inequality over $\tau\in (0,\varepsilon)$, $\epsilon\in(0,T^*-T)$ and divide the result by $\varepsilon$ to get
\begin{equation*}
    \begin{aligned}
    0&\le \intq (\tens^{N} - \tens^\ast(\vth^N, \dens )):(D\bu^{N} - \dens) \dq \le  \frac{1}{2}\left(\|\bu^N_0\|_2^2 - \frac{1}{\varepsilon}\int_0^{\varepsilon}\|\bu^N(T+\tau)\|_2^2\;\mathrm{d}\tau\right)  \\
    &+ \int_0^{T} \dual{\beh}{\bu^{N}}_V \dt - \int_0^{T}\intom \tens^{N}: \dens   +  \tens^\ast(\vth^N, \dens ):(D\bu^{N} - \dens ) \dq.
\end{aligned} 
\end{equation*}
Using the convergence results \eqref{Nconvergence3}--\eqref{Nconvergence5}, the convergence \eqref{aeNconvergence}, the assumption \eqref{tensor}, and the convergence of the initial conditions \eqref{icconvergences}, we can let $N \to \infty$ in the above inequality to obtain
\begin{equation*}
    \begin{aligned}
    0&\le  \lim_{N\to\infty} \intq (\tens^{N} - \tens^\ast(\vth^N, \dens )):(D\bu^{N} - \dens) \dq \le \frac{1}{2}\left(\|\bu_0\|_2^2 - \frac{1}{\varepsilon}\int_0^{\varepsilon}\|\bu(T+\tau)\|_2^2\;\mathrm{d}\tau\right)  \\
    &+ \int_0^{T} \dual{\beh}{\bu}_V \dt - \int_0^{T}\intom \tens: \dens   +  \tens^\ast(\vth, \dens ):(D\bu - \dens ) \dq\\
    &= \frac{1}{2}\left(\|\bu(T)\|_2^2 - \frac{1}{\varepsilon}\int_0^{\varepsilon}\|\bu(T+\tau)\|_2^2\;\mathrm{d}\tau\right) + (\tens- \tens^\ast(\vth, \dens )):(D\bu - \dens ) \dq,
\end{aligned} 
\end{equation*}
where for the last equality, we used the identity \eqref{limitequationE}. Letting $\varepsilon \to 0_+$ and using the continuity of $\bu \in C([0,T^*];L_{0,\diver}^2)$, we deduce
\begin{equation}\label{Minty}
    \begin{aligned}
    &0\le  \lim_{N\to\infty} \intq (\tens^{N} - \tens^\ast(\vth^N, \dens )):(D\bu^{N} - \dens) \dq \le \intq(\tens - \tens^\ast(\vth, \dens )):(D\bu - \dens ) \dq.
\end{aligned} 
\end{equation}
We now take $\dens :=D\bu \mp \varepsilon \wens$ for some $\wens \in L^p(Q;\rr^{2\times2}_{sym})$ and $\varepsilon\in(0,1)$. Hence, we have
\begin{align*}
    &0\leq \pm \intq (\tens - \tens^*(\vth, D\bu \mp \varepsilon \wens)): \wens\dq\stackrel{\varepsilon \to 0+}{\to}\pm \intq (\tens - \tens^*(\vth, D\bu )): \wens \dq, 
\end{align*} 
since $\tens^\ast$ is continuous. Consequently,  
\begin{align*}
    0 = \intq (\tens - \tens^*(\vth, D\bu )): \wens \dq
\end{align*} 
for all  $\wens \in L^p(Q;\rr^{2\times2}_{sym})$. Therefore, $\tens = \tens^*(\vth, D\bu )$ almost everywhere in $Q$.
We thus have $\bu \in L^p(0,T; V)$, such that  $\dert \bu \in L^{p'}(0,T; V^\ast)$ and 
\begin{equation}
  \begin{aligned}
       &\intt\dual{\dert \bu}{ \bw}_V\dt + \intq \tens^\ast(\vth, D\bu ):D\bw \dq = \intq (\bu \otimes \bu):\nabla\bw \dq   +   \intt\dual{\beh}{\bw}_V\dt, 
\end{aligned}   
\end{equation}
for all $\bw \in L^p(0,T; V)$. 
We have found $\bu$ fulfilling the \emph{Momentum equation}~\eqref{me}.

Next, we focus on the limit passage in the temperature equation. We need to strengthen the convergence result for $\tens^N$ and $D\bu^N$, and show that
\begin{align}\label{tensconvergence}
     &\tens^N:D\bu^{N} \rightharpoonup \tens:D\bu && \text{weakly in $L^{1}(Q)$}.
\end{align}
We set $\dens := D\bu$ in \eqref{Minty} and immediately obtain the convergence
\begin{align}\label{strongForMinty}
     &(\tens^{N} - \tens^*(\vth^N, D\bu )):(D\bu^{N} - D\bu ) \to 0 && \text{in $L^{1}(Q)$}.
\end{align}
Next, let $w \in L^\infty(Q)$ be arbitrary. Then, using \eqref{Nconvergence3}--\eqref{Nconvergence5}, \eqref{aeNconvergence}, the assumption \eqref{tensor}, and the dominated convergence theorem, we have
\begin{align}\label{pomoz}
     \tens^\ast(\vth^N, D\bu ) &\to \tens  && \text{strongly in $L^{p'}(Q; \mathbb{R}^{2\times 2})$},\\
     D\bu^N w &\rightharpoonup D\bu w  && \text{weakly in $L^{p}(Q; \mathbb{R}^{2\times 2})$}.\label{pomoz2}
\end{align}
Therefore
$$
\begin{aligned}
&\lim_{N\to \infty}   \intq \tens^N : D\bu^N w \dq = \lim_{N\to \infty} \intq \tens^N : (D\bu^N-D\bu) w + \tens^N : D\bu w \dq\\
&\overset{\eqref{Nconvergence5}}{=}  \lim_{N\to \infty} \intq (\tens^N- \tens^*(\vth^N, D\bu )): (D\bu^N-D\bu) w \dq\\
&\qquad  + \lim_{N\to \infty} \intq \tens^*(\vth^N, D\bu ) : (D\bu^N-D\bu) w \dq + \intq \tens : D\bu w \dq\\
&\overset{\eqref{strongForMinty}}{=} 0 + \lim_{N\to \infty} \intq (\tens^*(\vth^N, D\bu ) -\tens): (D\bu^N-D\bu) w \dq\\
&\qquad + \lim_{N\to \infty} \intq \tens: (D\bu^N-D\bu) w \dq +\intq \tens : D\bu w \dq \overset{\eqref{pomoz},\eqref{pomoz2}}{=} \intq \tens : D\bu w \dq 
\end{aligned}
$$
and hence \eqref{tensconvergence} holds.

%

Using the convergence result \eqref{tensconvergence}, we can now focus on the limiting procedure in the temperature equation \eqref{2gal}.
Let us now consider \eqref{2gal} with  $\psi \in C^{\infty}_0((-\infty,T)\times\Omega)$ and integrate it over $(0,T)$ to get
\begin{equation}\label{weakform2}
    \begin{aligned}
   &\intt\dual{\dert \vth^{N}}{\psi}_{W^{1,2}_0}\dt + \intq \kappa(\vth^{N}) \nabla\vth^{N} \cdot \nabla \psi  \dq =  \intq \tens^N:D\bu^{N}  \psi + \vth^{N} \bu^{N} \cdot \nabla \psi  \dq.
    \end{aligned} 
\end{equation}
Then, integrating by parts in the first term, we can use \eqref{Nconvergence1}--\eqref{strongthetaNconvergence}, \eqref{tensconvergence} and let $N \to \infty$ to conclude that
\begin{equation}\label{final2}
    \begin{aligned}
   &-\intq \vth \, \dert \psi \, \dq 
   + \intq \kappa(\vth) \nabla \vth \cdot \nabla \psi \, \dq 
   = \intq \tens : D\bu \, \psi + \vth \, \bu \cdot \nabla \psi \, \dq 
   + \intom \vth_0 \psi(0) \, \dx
    \end{aligned} 
\end{equation}
for all $\psi \in C^{\infty}_0((-\infty,T)\times\Omega)$. 

We have thus found $\vth \in L^\infty(0,T;L^1\om)\cap L^r(Q)$, for $r \in [1,2)$, such that 
\begin{equation}
\begin{aligned}
    &\vth - \htet \in L^s(0,T;W^{1,s}_0\om) \quad \text{for } s \in [1,4/3),\\
    &(\vth)^{\alpha} \in L^2(0,T;W^{1,2}\om) \quad \text{for } \alpha \in [0,1/2),
\end{aligned}\label{odkaz}
\end{equation}
and $\vth$ solves the \emph{Internal energy balance} \eqref{ieb}. 
Furthermore, the minimum principle $\vth^N \geq \mu$ for $N \in \mathbb{N}$, together with \eqref{aeNconvergence}, implies $\vth \geq \mu$ almost everywhere.

\subsection{Initial condition for temperature}
 Notice that \eqref{final2} already indicates a kind of weak attainment of the initial condition (without specifying it rigorously). Here, our goal is different; we show the strong attainment in $L^1(\Omega)$, but only for the essential limit as $t \to 0_+$. In the next section, when we show that the temperature is continuous, we prove the attainment of the initial condition for the classical limit $t \to 0_+$. The goal here is to show the existence of a set $S \subset (0,T)$ such that $(0,T) \setminus S$ is of measure zero and
\begin{align}\label{thetaticae}
     &\vth(t)\to \vth_0 &&\text{in $L^1(\Omega)$ as $S \ni t\to0$.}
 \end{align}
 
From \eqref{strongthetaNconvergence}, we know that there exists a subsequence of $\vth^N$ and a set $S\subset(0,T) $ with $(0,T)\setminus S$ of measure zero such that, for all $t\in S$, it holds 
\begin{align}\label{Strik}
    &\vth^N(t)\to\vth(t)&&\text{in $L^s(\Omega)$, $s<2$.}
\end{align}
Let us now realize that by \eqref{LinftyL1} we have $\sqrt{\theta(t)}\in L^2(\Omega)$ for a.e. $t\in(0,T)$ and show that
\begin{align}\label{L2sqrtconvineq}
    &\liminf_{S\ni t\to0+} \intom \sqrt{\vth(t)}\varphi \dx \geq \intom \sqrt{\vth_0}\varphi \dx &&\text{ for all $\varphi\in L^2(\Omega)$, $\varphi \geq 0$}.
\end{align}
For $0 \leq \varphi \in C^\infty_0(\Omega)$, we set $\psi := \frac{\varphi}{\sqrt{\vth^N}}$ in equation \eqref{2gal}, and integrate over the time interval $(0,t)$, with $t \in S$ arbitrary. Note that such $\psi$ is an admissible test function for \eqref{2gal}, since $\frac{\varphi}{\sqrt{\vth^N}}(t) \in W^{1,2}_0(\Omega)$ for almost all $t \in (0,T)$. Additionally, $\frac{\varphi}{\sqrt{\vth^N}} \in L^\infty(Q)$ due to the minimum principle. We obtain
\begin{equation*}
    \begin{aligned}
   &\int_0^t  \dert \intom 2\sqrt{\vth^N}\varphi\dx \dtau + \int_0^t \intom \kappa(\vth^{N}) \nabla\vth^{N} \cdot \nabla \frac{\varphi}{\sqrt{\vth^N}}  \dx \dtau= &&\\
    & = \int_0^t \intom \tens^{N}:D\bu^{N} \frac{\varphi}{\sqrt{\vth^N}}  \dx \dtau + \int_0^t \intom \vth^{N} \bu^{N} \cdot \nabla \frac{\varphi}{\sqrt{\vth^N}} \dx \dtau,
    \end{aligned} 
\end{equation*}
using integration by parts and the fact that $\diver \bu^N=0$, we can rewrite it as
\begin{equation}\label{sqrtgal}
    \begin{aligned}
   & 2\intom \sqrt{\vth^N(t)}\varphi-\sqrt{\vth^N_0}\varphi\dx + \int_0^t\intom \left(\frac{\kappa(\vth^{N})}{\sqrt{\vth^N}}{\nabla\vth^{N}}-2{\sqrt{\vth^N}} \bu^N\right) \cdot  \nabla \varphi \dx \dtau= \\
    &= \int_0^t \intom \kappa(\vth^{N}) \varphi\frac{|\nabla\vth^{N}|^2}{2(\vth^N)^{3/2}}  \dx \dtau + \int_0^t \intom \tens^{N}:D\bu^{N} \frac{\varphi}{\sqrt{\vth^N}}  \dx \dtau \geq 0.
    \end{aligned} 
\end{equation}
We can now let $N\to \infty$ in the inequality \eqref{sqrtgal}. Using convergence results \eqref{Strik}, \eqref{icconvergences}, \eqref{aeNconvergence}, \eqref{Nconvergence7}, and \eqref{Nconvergence4} yields 
\begin{equation}\label{sqrtgallim}
    \begin{aligned}
   & 2\intom \sqrt{\vth(t)}\varphi-\sqrt{\vth_0}\varphi\dx + \int_0^t\intom \left(\frac{\kappa(\vth)}{\sqrt{\vth}} \nabla\vth-2{\sqrt{\vth}} \bu\right) \cdot  \nabla \varphi \dx \dtau \geq 0.
    \end{aligned} 
\end{equation}
We can pass to the essential limit inferior $S\ni t \to 0_+$ in \eqref{sqrtgallim}. The time integral vanishes, and we deduce
\begin{align}\label{weakM}
    &\liminf_{S\ni t\to 0_+} \intom \sqrt{\vth(t)}\varphi \dx \geq \intom \sqrt{\vth_0}\varphi \dx &&\text{ for all $\varphi\in C^\infty_0(\Omega)$, $\varphi\geq 0$}.
\end{align}
Since $\sqrt{\vth}\in L^\infty(0,T;L^2\om)$  and $C^\infty_0(\Omega)$ is dense in $L^2(\Omega)$, we can extend the result to all nonnegative $\varphi \in L^2(\Omega)$. Hence~\eqref{L2sqrtconvineq} holds. 
 
Next, we want to show 
\begin{align}\label{dualphiconvergence}
     &\lim_{S\ni t\to0+}\intom \vth(t)\varphi \dx = \intom \vth_0 \varphi &&\text{for $\varphi\in C^\infty_0(\Omega)$, $0\leq \varphi\leq 1$}.
\end{align}
We set $\psi:=\varphi$ in \eqref{2gal}, where   $\varphi\in C^\infty_0(\Omega)$, $0\leq \varphi\leq 1$ and integrate the result over time interval $(0,t)$, where $t\in S$ is arbitrary.  We have
\begin{equation}\label{Ncontinuity}
    \begin{aligned}
    &\intom \vth^{N}(t){\varphi} \dx =  \intom \vth^{N}_0{\varphi}\dx + \int_0^t\intom \tens^{N}:D\bu^{N} \varphi + \vth^{N} \bu^{N} \cdot \nabla \varphi - \kappa(\vth^{N}) \nabla\vth^{N} \cdot \nabla \varphi \dx\dtau.
    \end{aligned} 
\end{equation}
We now want to pass to the limit $N\to \infty$ in \eqref{Ncontinuity}. The right-hand side converges by \eqref{aeNconvergence}, \eqref{Nconvergence7}, \eqref{icconvergences}, \eqref{tensconvergence}, \eqref{Nconvergence4} and \eqref{Nconvergence6}. The left-hand side converges by \eqref{Strik}.
We thus have 
\begin{equation}\label{cont}
    \begin{aligned}
    \intom \vth(t){\varphi} \dx =  \intom \vth_0{\varphi}\dx + \int_0^t\intom \tens:D\bu \varphi + \vth \bu \cdot \nabla \varphi - \kappa(\vth) \nabla\vth \cdot \nabla \varphi \dx\dtau 
    \end{aligned} 
\end{equation}
We now let $S \ni t\to0$ in \eqref{cont} %
and obtain the desired result $\eqref{dualphiconvergence}$. 

Finally, we show that \eqref{L2sqrtconvineq} and \eqref{dualphiconvergence} imply the attainment of the initial condition locally in $\Omega$. In fact, for $\varphi\in C^\infty_0(\Omega)$, $0\leq \varphi\leq 1$, we can estimate
\begin{align*}
     0&\leq \lim_{S\ni t\to0+}\intom \left(\sqrt{\vth(t)} - \sqrt{\vth_0}\right)^2{\varphi}\dx \leq \lim_{S \ni t\to0+} \intom \vth(t)\varphi + \vth_0{\varphi}\dx - \liminf_{S \ni t\to0+} 2\intom \sqrt{\vth(t)}\sqrt{\vth_0}\varphi\dx\\
     &\stackrel{\eqref{dualphiconvergence},\eqref{L2sqrtconvineq}}{\leq} 0, 
\end{align*}
where we can use the convergence result \eqref{L2sqrtconvineq}, since $\sqrt{\vth_0}\varphi$ belongs to $L^2(\Omega)$.  
Let $K \subset \Omega$ be an arbitrary compact set. There exists a function $\varphi \in C^\infty_0(\Omega)$, with $0 \le \varphi \le 1$, such that $\varphi \equiv 1$ on $K$.
Hence, the inequality above therefore yields that, as $S \ni t \to 0$,
\begin{align*}
\vth^{1/2}(t) \to \vth^{1/2}_0
\qquad \text{in } L^{2}(K), \text{ for any compact } K \subset \Omega.
\end{align*}
Due to this strong convergence, for any $\{t_k\}_k\subset S$ with $t_k\to0$ as $k\to+\infty$, the sequence $\{\vth(t_k)\}_k$ is uniformly equiintegrable and one can extract $\{k_n\}_n$ such that $\vth(t_{k_n}) \to \vth_0$ almost everywhere in $K$ as $n \to +\infty$. Consequently, for any $\beta \in (0,1]$,
\begin{align}\label{strong-konv-on-K}
\vth^{\beta}(t) \to \vth^{\beta}_0
\qquad \text{in } L^{\frac{1}{\beta}}(K), \text{ for any compact } K \subset \Omega,
\end{align}
as $S \ni t \to 0$.
This and the uniform bound \eqref{LinftyL1} then turn into
\begin{align}\label{sqrtconvergenceth0}
    &\vth^{\beta}(t) \rightharpoonup \vth^{\beta}_0 &&\text{in $L^{\frac{1}{\beta}}(\Omega)$ as $S \ni t\to0$.}
\end{align}

Lastly, it is enough to show that
\begin{align}\label{2sqrttwiceconvergence}
    &\limsup_{S \ni t\to0} \int_\Omega \vth(t) \, \mathrm{d}x \leq \int_\Omega \vth_0 \, \mathrm{d}x,
\end{align}
which, combined with \eqref{sqrtconvergenceth0}, leads to the desired result
\begin{align}\label{sqrtconvergenceth0A}
    &\vth^{\beta}(t) \to \vth^{\beta}_0 &&\text{strongly in $L^{\frac{1}{\beta}}(\Omega)$ as $S \ni t\to0$,}
\end{align}
valid for all $\beta \in (0,1]$. We can proceed here as when deriving \eqref{strong-konv-on-K}. We first consider $\beta=1/2$ and then the general $\beta$.
 
Hence, we focus on  \eqref{2sqrttwiceconvergence}. We set $\psi:=1-\left(\frac{\htet}{\vth^N}\right)^{\frac{\alpha+1}{2}}$ in \eqref{2gal} and integrate over time interval $(0,t)$, where $t\in S$ is arbitrary. Note, that $\psi$ is an admissible test function since it is zero at the boundary, $\psi(t) \in L^\infty(\Omega)$ for a.a. $t\in S$, and
$$
\nabla \psi =-\frac{\alpha+1}{2}\left(\frac{\htet}{\vth^N}\right)^{\frac{\alpha+1}{2}}\left(\frac{\nabla \htet}{\htet} -\frac{\nabla \vth^N}{\vth^N} \right)\in L^2(\Omega) \text{ almost everywhere in $S$,}
$$
which follows from \eqref{odkaz} and \eqref{conditionstheta}  for all $\alpha>0$. We obtain
\begin{equation*}
    \begin{aligned}
   &\intom \left(\vth^N - \frac{2(\vth^N)^{\frac{1-\alpha}{2}}(\htet)^{\frac{\alpha+1}{2}}}{1-\alpha}\right)(t) \dx + \intu{t} \intom \kappa(\vth^{N}) \nabla\vth^{N} \cdot \nabla \psi \dx\dtau = \\
    & = \intu{t} \intom \tens^N:D\bu^{N} \psi  +  \vth^{N} \bu^{N} \cdot \nabla \psi  \dx \dtau + \intom \left(\vth^N_0 - \frac{2(\vth^N_0)^{\frac{1-\alpha}{2}}(\htet)^{\frac{\alpha+1}{2}}}{1-\alpha}\right) \dx. 
    \end{aligned} 
\end{equation*}
Neglecting the nonnegative  term 
$$\intu{t}\intom \frac{\kappa(\vth^{N})(\alpha+1)(\htet)^{\frac{\alpha+1}{2}}}{2(\vth^N)^{\frac{\alpha+3}{2}}}|\nabla\vth^{N}|^2 \dx\dtau,$$
and using the fact that $\psi$ is bounded, we get 
\begin{equation}\label{inequalitySN}
    \begin{aligned}
   & \intom \left(\vth^N(t) - \frac{2(\vth^N(t))^{\frac{1-\alpha}{2}}(\htet)^{\frac{\alpha+1}{2}}}{1-\alpha}\right) \dx  \leq   \intom \left(\vth^N_0 - \frac{2(\vth^N_0)^{\frac{1-\alpha}{2}}(\htet)^{\frac{\alpha+1}{2}}}{1-\alpha}\right) \dx + C\intu{t} \intom \tens^N:D\bu^{N} \dx\dtau \\ 
   &\qquad +\frac{\alpha+1}{1-\alpha}\intu{t} \intom  \kappa(\vth^{N})(\htet)^{\frac{\alpha-1}{2}} \nabla \htet\cdot  \nabla (\vth^{N})^{\frac{1-\alpha}{2}}  \dx\dtau  + \intu{t} \intom \vth^N \bu^N \cdot \nabla \left(\frac{\htet}{\vth^N}\right)^{\frac{\alpha+1}{2}} \dx \dtau. 
    \end{aligned} 
\end{equation}
We want to pass to the limit $N \to \infty$ in \eqref{inequalitySN}. The term on the left-hand side converges by \eqref{Strik}. The first term on the right-hand side converges by \eqref{icconvergences}. The convergence of the second term on the right-hand side holds due to \eqref{tensconvergence}. For the third term, it is sufficient to use the weak convergence result \eqref{Nconvergence8}, the assumptions \eqref{kap}, \eqref{conditionstheta}, and the pointwise convergence \eqref{aeNconvergence}.
%
We use integration by parts, the fact that $\diver \bu^N = 0$, and rewrite the fourth term on the right-hand side of \eqref{inequalitySN} as follows
$$
\begin{aligned}
&\intu{t} \intom \vth^N \bu^N \cdot \nabla \left(\frac{\htet}{\vth^N}\right)^{\frac{\alpha+1}{2}} \dx \dtau
= \frac{\alpha+1}{2}\intu{t}\intom \htet \left(\frac{\htet}{\vth^N}\right)^{\frac{\alpha-3}{2}} \bu^N \cdot \nabla \left(\frac{\htet}{\vth^N}\right) \dx \dtau\\
&\quad = \frac{\alpha+1}{\alpha-1}\intu{t}\intom \htet \bu^N \cdot \nabla \left(\frac{\htet}{\vth^N}\right)^{\frac{\alpha-1}{2}} \dx \dtau
= \frac{\alpha+1}{1-\alpha}\intu{t}\intom (\vth^N)^{\frac{1-\alpha}{2}} \bu^N \cdot \frac{\nabla \htet}{(\htet)^{\frac{1-\alpha}{2}}} \dx \dtau.
\end{aligned}
$$
The resulting integral can be easily identified in the limit due to the assumptions \eqref{conditionstheta} and \eqref{muconditions}, the convergence results \eqref{strongthetaNconvergence} and \eqref{Nconvergence4}, and the facts that $\alpha \in (0,1)$ and $p \ge 2$.  
We thus obtain
\begin{equation}\label{inequalityS}
\begin{aligned}
& \intom \left(\vth(t) - \frac{2(\vth(t))^{\frac{1-\alpha}{2}}(\htet)^{\frac{\alpha+1}{2}}}{1-\alpha}\right) \dx  
\leq \intom \left(\vth_0 - \frac{2(\vth_0)^{\frac{1-\alpha}{2}}(\htet)^{\frac{\alpha+1}{2}}}{1-\alpha}\right) \dx 
+ C\intu{t} \intom \tens : D\bu \dx \dtau \\ 
&\qquad + \frac{\alpha+1}{1-\alpha}\intu{t} \intom \kappa(\vth) (\htet)^{\frac{\alpha-1}{2}} \nabla \htet \cdot \nabla (\vth)^{\frac{1-\alpha}{2}} \dx \dtau 
+ \intu{t} \intom \vth \bu \cdot \nabla \left(\frac{\htet}{\vth}\right)^{\frac{\alpha+1}{2}} \dx \dtau.
\end{aligned}
\end{equation}
Let us take the limit $S \ni t\to0$ in \eqref{inequalityS}. Since all the terms on the right-hand side of \eqref{inequalityS} except for the first one converge to zero, we obtain 
\begin{align*}
     &\limsup_{S \ni t\to0}\intom \left(\vth(t) - \frac{2(\vth(t))^{\frac{1-\alpha}{2}}(\htet)^{\frac{\alpha+1}{2}}}{1-\alpha}\right) \dx  \leq   \intom \left(\vth_0 - \frac{2(\vth_0)^{\frac{1-\alpha}{2}}(\htet)^{\frac{\alpha+1}{2}}}{1-\alpha}\right) \dx.
\end{align*}
Using \eqref{sqrtconvergenceth0} with $\beta := \frac{1-\alpha}{2}$, it follows from the above inequality that 
\begin{align*}
     &\limsup_{S \ni t\to0}\intom \vth(t) \dx  \leq   \intom \vth_0 \dx,
\end{align*}
which is precisely the relation \eqref{2sqrttwiceconvergence} we aimed to prove. Consequently, we also deduce the strong attainment of the initial condition \eqref{sqrtconvergenceth0A}.

\subsection{Entropy equality} 
Proving that the solution found in the previous sections fulfills the equality \eqref{ee} can follow the approach of \cite{Abbatiello2022}, i.e., we set  $\psi := \frac{\varphi}{\vth^N}$, where $\varphi \in C^\infty_0((-\infty,T)\times \Omega)$, in \eqref{2gal} and integrate the resulting expression over $(0,T)$ to get
\begin{equation}\label{entropygal}
    \begin{aligned}
     \intq \kappa(\vth^N) \nabla\eta^N \cdot \nabla \varphi 
 - \eta^N \dert \varphi  - \eta^N \bu^N \cdot \nabla \varphi\dq 
    = \intq \frac{\tens^N:D\bu^N \varphi}{\vth^N}  + \kappa(\vth^N)
    \frac{|\nabla\vth^N|^2}{(\vth^N)^2} \varphi \dq + \intom \eta^N_0 \varphi(0) \dx,
    \end{aligned}
\end{equation}
where $\eta^N_0 := \log \vth^N_0$ and $\eta^N := \log \vth^N$. Next, the convergence as $N \to \infty$ is shown. In what follows, we only highlight the key convergence results that are used in the next section and refer the reader to the original paper for a detailed proof. Note that to prove the entropy equality, we need the convergence of $\vth$ to the initial condition in the sense of \eqref{thetaticae}.

By \eqref{strongthetaNconvergence}, for $N\to\infty$, we have 
\begin{align}\label{etastrong}
    &\eta^{N} \to \eta && \text{in $L^{q}(Q)$ for $q\in \left[1,\infty\right)$},
\end{align}
where $\eta := \log \vth$.
Moreover,
$\nabla\eta^{N} = \frac{\nabla \vth^N}{\vth^N}$ is bounded in $L^2(Q)$ by \eqref{estimatecombofinal}. We can thus find a converging subsequence 
\begin{align}\label{etagrad}
    &\nabla\eta^{N} \rightharpoonup \nabla\eta && \text{in $L^{2}(Q)$},
\end{align}
and so $\eta \in L^2(0,T;W^{1,2}\om)\cap L^q(Q)$ for $q\in \left[1,\infty\right)$.  
Following \cite[Subsection 3.(b)(i)]{Abbatiello2022}, 
we have 
\begin{align}
    &\nabla \vth^N\to \nabla\vth &&\text{ in }L^1(Q)\\
    &\nabla \vth^N \to \nabla \vth &&\text{ a.e. in }Q.\label{gradthetpoint}
\end{align}

Moreover, from  \cite[Subsection 3.(b)(ii)]{Abbatiello2022}, 
we know
\begin{align}\label{NstrongL1entropy}
    &\kappa(\vth^N)\frac{|\nabla\vth^N|^2}{(\vth^N)^2}\to {\kappa(\vth)}\frac{|\nabla\vth|^2}{\vth^2} &&\text{in $L^1(Q)$.}
\end{align}
The convergence results \eqref{etastrong}--\eqref{NstrongL1entropy} are enough to pass to the limit in \eqref{entropygal}.
Consequently, $\eta \in L^2(0,T;W^{1,2}\om)\cap L^q(Q), q\in \left[1,\infty\right)$ satisfies \emph{Entropy equation} \eqref{ee}.

\subsection{Continuity}

This section is the final part of the proof. It relies on the approach developed in \cite{Abbatiello2022}, but here we are able to show a stronger result, namely the continuity with respect to time. Our goal is to show that the sequence $\vth^N$ is Cauchy in $C([0,T]; L^1(\Omega))$, and therefore its limit shares this property. The proof is divided into two parts. First, we establish this property for the truncated functions, and second, we prove a kind of uniform equiintegrability.

Recall the definitions of truncations \eqref{tk} and \eqref{propertiestke}. Our first goal is to show that $\{\mathcal{T}_{m,\delta}(\vth^N)\}_{N}$ is Cauchy in $C([0,T];L^2\om)$ for any fixed $m>\norm{\htet}{\infty}{ }$ and~$\delta\in (0,m/2)$. This means that for any given $\alpha>0$, there exists $N_0\in\mathbb{N}$ such that for all $N,M > N_0$ it holds
\begin{equation}\label{goal}
\sup_{t\in [0,T]} \|\mathcal{T}_{m,\delta}(\vth^N(t))-\mathcal{T}_{m,\delta}(\vth^M)(t)\|_2^2 \le \alpha. 
\end{equation}
We set $\psi^N:=\mathcal{T}'_{m,\delta}(\vth^N)(\mathcal{T}_{m,\delta}(\vth^N) - \mathcal{T}_{m,\delta}(\vth^M))$, $\psi^M:=\mathcal{T}'_{m,\delta}(\vth^M)(\mathcal{T}_{m,\delta}(\vth^N) - \mathcal{T}_{m,\delta}(\vth^M))$ in \eqref{2gal} for $N$, $M\in \nn$, respectively. Subtracting one equation from the other and integrating the result over time interval $(0,t)$ with arbitrary $t\in [0,T]$, we obtain \begin{equation}\label{tktkesti}
    \begin{aligned}
    &\frac{1}{2}\norm{\mathcal{T}_{m,\delta}(\vth^N(t)) - \mathcal{T}_{m,\delta}(\vth^M(t))}{2}{2}= \frac{1}{2}\norm{\mathcal{T}_{m,\delta}(\vth^N_0) - \mathcal{T}_{m,\delta}(\vth^M_0)}{2}{2}\\        
    & \quad +\int_0^t \intom (\mathcal{T}'_{m,\delta}(\vth^N)\tens^{N}:D\bu^{N}-\mathcal{T}'_{m,\delta}(\vth^M)\tens^{M}:D\bu^{M})(\mathcal{T}_{m,\delta}(\vth^N) - \mathcal{T}_{m,\delta}(\vth^M))\dx \dtau\\
    &\quad+ \int_0^t\intom  \big(\mathcal{T}_{m,\delta}(\vth^N) \bu^N - \mathcal{T}_{m,\delta}(\vth^M)\bu^M \big)\cdot\nabla(\mathcal{T}_{m,\delta}(\vth^N) -  \mathcal{T}_{m,\delta}(\vth^M)) \dx \dtau\\
     &\quad  - \int_0^t \intom (\kappa(\vth^{N}) \nabla \mathcal{T}_{m,\delta}(\vth^N) - \kappa(\vth^{M}) \nabla \mathcal{T}_{m,\delta}(\vth^M))\cdot \nabla(\mathcal{T}_{m,\delta}(\vth^N) - \mathcal{T}_{m,\delta}(\vth^M))\dx \dtau\\ 
     &\quad  - \int_0^t \intom \bigg(\kappa(\vth^{N}) \mathcal{T}''_{m,\delta}(\vth^N)|\nabla\vth^N|^2 - \kappa(\vth^{M}) \mathcal{T}''_{m,\delta}(\vth^M)|\nabla\vth^M|^2\bigg)(\mathcal{T}_{m,\delta}(\vth^N) - \mathcal{T}_{m,\delta}(\vth^M))\dx \dtau. 
    \end{aligned} 
\end{equation}
Using the simple estimate  $|\mathcal{T}'_{m,\delta}(\vth^N)\vth^{N}|\le m+\delta$ and the assumption \eqref{kap}, we deduce
\begin{equation}\label{cauhytk}
    \begin{aligned}
    &\sup_{t\in[0,T]}\norm{\mathcal{T}_{m,\delta}(\vth^N(t)) - \mathcal{T}_{m,\delta}(\vth^M(t))}{2}{2} \leq 2m\norm{\mathcal{T}_{m,\delta}(\vth^N_0) - \mathcal{T}_{m,\delta}(\vth^M_0)}{1}{ }\\     
    & +2\int_0^T \intom \left(\left|{\tens^{N}:D\bu^{N}}\right|+\left|\tens^{M}:D\bu^{M}\right|\right)\left|\mathcal{T}_{m,\delta}(\vth^N) - \mathcal{T}_{m,\delta}(\vth^M)\right|\dx \dt\\
    &+ 4Cm\int_0^T\intom  \left(|\bu^N|+ |\bu^M| + |\nabla \mathcal{T}_{m,\delta}(\vth^N)|+ | \nabla \mathcal{T}_{m,\delta}(\vth^M)|\right)\left| \nabla(\mathcal{T}_{m,\delta}(\vth^N) - \mathcal{T}_{m,\delta}(\vth^M))\right| \dx \dt.
    \end{aligned} 
\end{equation}
Next, we estimate the right-hand side of \eqref{cauhytk}. By \eqref{icconvergences}, we have
\begin{equation}\label{cauchy1}
    2m\norm{\mathcal{T}_{m,\delta}(\vth^N_0) - \mathcal{T}_{m,\delta}(\vth^M_0)}{1}{ }<\frac{\alpha}{3}
\end{equation}
for $M$, $N$ large enough.
To estimate the second term, we can split 
 \begin{align*}
     \left|\mathcal{T}_{m,\delta}(\vth^N) - \mathcal{T}_{m,\delta}(\vth^M)\right| = \left|\mathcal{T}_{m,\delta}(\vth^N) - \mathcal{T}_{m,\delta}(\vth)\right| + \left|\mathcal{T}_{m,\delta}(\vth) - \mathcal{T}_{m,\delta}(\vth^M)\right|.
 \end{align*}
We know that for all $N\in\nn$, the function ${\mathcal{T}_{m,\delta}(\vth) - \mathcal{T}_{m,\delta}(\vth^N)}$ is bounded in $L^\infty(Q)$ by \eqref{propertiestke}, and that it converges almost everywhere to zero by the continuity of $\mathcal{T}_{m,\delta}$ and the convergence \eqref{aeNconvergence}.
Furthermore, by~\eqref{tensconvergence},  $\tens^{N}:D\bu^{N}$ is uniformly integrable and bounded in $L^1(Q)$ for all $N\in\nn$. Using Egorov's theorem, we can conclude that 
\begin{align}
    2\int_0^T \intom \left(\left|{\tens^{N}:D\bu^{N}}\right|+\left|\tens^{M}:D\bu^{M}\right|\right)\left|\mathcal{T}_{m,\delta}(\vth^N) - \mathcal{T}_{m,\delta}(\vth^M)\right|\dx \dt \leq \frac{\alpha}{3}
\end{align}
for $M,N$ large enough. Since $\kappa(\vth^{N}) \mathcal{T}''_{m,\delta}(\vth^N)|\nabla\vth^N|^2\to \kappa(\vth) \mathcal{T}''_{m,\delta}(\vth)|\nabla\vth|^2$ in $L^1(Q)$ as $N\to+\infty$ by \eqref{aeNconvergence}, \eqref{gradthetpoint} and \eqref{NstrongL1entropy}, we can also deal in the same manner with the last term on the right hand side to get that for large $N,M$ it can be estimated by $\alpha/3$. 
For the second and third terms on the right-hand side of \eqref{tktkesti}, we use that
$$|\bu^N|+ | \bu^M| + |\nabla \mathcal{T}_{m,\delta}(\vth^N)|+ | \nabla \mathcal{T}_{m,\delta}(\vth^M)|$$ is bounded in $L^2(Q)$, and 
$$\left| \nabla(\mathcal{T}_{m,\delta}(\vth^N) - \mathcal{T}_{m,\delta}(\vth^M))\right|\leq \left| \nabla(\mathcal{T}_{m,\delta}(\vth^N) - \mathcal{T}_{m,\delta}(\vth))\right|+\left| \nabla(\mathcal{T}_{m,\delta}(\vth) - \mathcal{T}_{m,\delta}(\vth^M))\right|,$$
where each term converges strongly to zero in $L^2(Q)$ by \eqref{NstrongL1entropy} and \eqref{aeNconvergence}.
For $M,N$ large enough, we can conclude that
\begin{align}\label{cauchy3}
    2\int_0^T\intom  \left(2m|\bu^N|+ 2m|\bu^M| + \overline{\kappa}|\nabla \mathcal{T}_{m,\delta}(\vth^N)|+ \overline{\kappa}| \nabla \mathcal{T}_{m,\delta}(\vth^M)|\right)\left| \nabla(\mathcal{T}_{m,\delta}(\vth^N) - \mathcal{T}_{m,\delta}(\vth^M))\right| \dx \dt<\frac{\alpha}{3}.
\end{align}
The sequence $\{\mathcal{T}_{m,\delta}(\vth^N)\}_{N}$ is therefore Cauchy in $C([0,T];L^2\om)$, since by \eqref{cauhytk}--\eqref{cauchy3} we have
\begin{align}\label{cauchy}
    \sup_{t\in[0,T]}\norm{\mathcal{T}_{m,\delta}(\vth^N(t)) - \mathcal{T}_{m,\delta}(\vth^M(t))}{2}{2} < \alpha,
\end{align}
for all $M,N$ sufficiently large.

Next, we recall $m>\norm{\htet}{\infty}{ }$, consider equation \eqref{2gal} with $\psi := \frac{(\vth^N-m)_+}{\vth^N}$, and integrate over time interval~$(0,t)$. Note that $\psi$ has zero trace due to the choice of $m$. This gives us 
\begin{equation}\label{continuitymid1}
    \begin{aligned}
         &\int_0^t \intom  \dert \vth^{N}\left( 1 - \frac{m}{\vth^N}\right) \chi_{\{\vth^N>m\}}\dx \dtau + \int_0^t \intom \frac{m\kappa(\vth^{N}) }{(\vth^N)^2}|\nabla\vth^{N}|^2\chi_{\{\vth^N>m\}}\dx \dtau = &&\\
        & = \int_0^t \intom \tens^{N}:D\bu^{N}\left(1-\frac{m}{\vth^N}\right)\chi_{\{\vth^N>m\}}\dx \dtau- \int_0^t\intom \bu^{N} \cdot \nabla\vth^N  \frac{(\vth^N-m)_+}{\vth^N} \dx \dtau. 
    \end{aligned}
\end{equation}
We can rewrite
$$\dert \vth^{N}\left( 1 - \frac{m}{\vth^N}\right) \chi_{\{\vth^N>m\}}= \dert( (\vth^{N} - m\log \vth^N -m + m\log m)_+ ).$$
Furthermore, the second term on the left-hand side of \eqref{continuitymid1} is positive, so we can omit it and replace the equality sign with an inequality. Since $\diver \vth^N =0$, we can use integration by parts to get rid  of the last term on the right-hand side. 
Hence, using $$0\leq\left(1-\frac{m}{\vth^N}\right)\chi_{\{\vth^N>m\}}\leq \chi_{\{\vth^N>m\}},$$
and the properties \eqref{tensor}, we can rewrite \eqref{continuitymid1} as
\begin{align*}
    \int_0^t \dert \norm{(\vth^{N} - m\log \vth^N - m + m\log m)_+}{1}{ } \dtau \leq \int_0^t \intom \tens^{N}:D\bu^{N} \chi_{\{\vth^N>m\}}\dx \dtau.
\end{align*}
Since $t \in (0,T)$ was chosen arbitrarily, we have
\begin{equation}\label{continuitymid2}
\begin{aligned}
    &\sup_{t\in[0,T]} \norm{(\vth^{N}(t) - m\log \vth^N(t) - m + m\log m)_+}{1}{ }  \leq \norm{(\vth^{N}_0 - m\log \vth^N_0 - m + m\log m)_+}{1}{ }\\
    &+\intt \intom \tens^{N}:D\bu^{N} \chi_{\{\vth^N>m\}}\dx \dt \leq \norm{\vth^{N}_0\chi_{\{\vth^N>m\}}}{1}{ } +\intt \intom \tens^{N}:D\bu^{N} \chi_{\{\vth^N>m\}}\dx \dt.
\end{aligned}
\end{equation}
Since $\tens^{N}:D\bu^{N}$ converges weakly to $\tens:D\bu$ in $L^1(Q)$ by \eqref{tensconvergence}, it forms a uniformly integrable sequence of functions. Additionally, due to the convergence \eqref{strongthetaNconvergence}, we have $|\{\vth^N>m\}|\le C/m$. This estimate on level sets of $\vth$, the equiintegrability of $\tens^{N}:D\bu^{N}$, and \eqref{icconvergences} imply that for any $\alpha>0$, there exists $m_0>0$ such that
\begin{align*}
   \norm{\vth^{N}_0\chi_{\{\vth^N_0>m\}}}{1}{ } +\intt \intom \tens^{N}:D\bu^{N} \chi_{\{\vth^N>m\}}\dx \dt < \frac{\alpha}{2}
\end{align*}
for all $m>m_0$.
For any $S$ positive, we can estimate $\frac{(S-m)_+}{2} \leq (S - m\log S - m + m\log m)_+$. The estimate \eqref{continuitymid2} therefore yields that for any $\alpha>0$ there exists $m_0>0$ such that 
\begin{align}\label{continuitymiddone}
   \sup_{t\in[0,T]} \norm{(\vth^{N}(t) - m)_+}{1}{ } < \alpha
\end{align}
for all $m>m_0$.

Finally, we show that the sequence $\{\vth^N\}$ is a Cauchy sequence in $C([0,T];L^1\om)$. Note that each element of the sequence belongs to $C([0,T];L^2\om)\subset C([0,T];L^1\om)$ by \eqref{thetaNcontinuity}. We recall that we consider $m>\norm{\htet}{\infty}{ }$ and~$\delta<m/2$.
For any $t\in(0,T)$ we can estimate
\begin{equation}\label{continuityfinal}
    \begin{aligned}
        \norm{\vth^{N}(t) - \vth^M(t)}{1}{ } &\leq \norm{\vth^{N}(t) - \mathcal{T}_{m,\delta}(\vth^N)(t) }{1}{ } + \norm{\vth^{M}(t) - \mathcal{T}_{m,\delta}(\vth^M)(t)}{1}{ }+ \norm{\mathcal{T}_{m,\delta}(\vth^N)(t) - \mathcal{T}_{m,\delta}(\vth^M)(t)}{1}{ }\\
        &\leq \norm{(\vth^{N}(t) - \frac{m}{2})_+}{1}{ } + \norm{(\vth^{M}(t) - \frac{m}{2})_+}{1}{ }  + \norm{\mathcal{T}_{m,\delta}(\vth^N)(t) - \mathcal{T}_{m,\delta}(\vth^M)(t)}{1}{ }
    \end{aligned}
\end{equation}
Hence, for a given $\varepsilon>0$, we can use \eqref{continuitymiddone} to find $m>2$ sufficiently large such that for all $M,N$
\begin{align}
     &\sup_{t\in[0,T]}\norm{(\vth^{N}(t) - \frac{m}{2})_+}{1}{ }<\frac{\varepsilon}{3}, &&\sup_{t\in[0,T]}\norm{(\vth^{M}(t) - \frac{m}{2})_+}{1}{ } <\frac{\varepsilon}{3}.
\end{align}
Finally, for this choice of $m$, we set $\delta = 1$ and find $N_0$ such that for all $M, N > N_0$ we have  
\begin{align*}
    \|\mathcal{T}_{m,\delta}(\vth^N)(t) - \mathcal{T}_{m,\delta}(\vth^M)(t)\|_{1} < \frac{\varepsilon}{3},
\end{align*}
by virtue of \eqref{cauchy}. We conclude that for any $\varepsilon>0$ there exists $N_0$ such that for all $M,N > N_0$ it holds
\begin{equation}
   \sup_{t \in (0,T)} \norm{\vth^{N}(t) - \vth^M(t)}{1}{} < \varepsilon.
\end{equation}
Hence, the sequence $\{\vth^N\}_{N=0}^\infty$ is Cauchy in $C([0,T];L^1\om)$, and its limit $\vth$ indeed belongs to $C([0,T];L^1\om)$.

\section{Proof of the role of the entropy equality - Theorem~\ref{maintheoremB}}\label{sec.4}
To establish the continuity properties of $\vth$, we first show that the entropy $\eta$, satisfying the assumptions of Theorem~\ref{maintheoremB}, belongs to $C([0,T];L^q)$ for all $q \in [1,\infty)$. The key advantage of the entropy identity~\eqref{eeB} is that it admits a straightforward renormalization; we provide full details for the sake of completeness. 

Using the density argument together with assumptions \eqref{vb1}--\eqref{vb4}, we observe that \eqref{eeB} extends to test functions of the form $\varphi(t,x):= \psi(t)\zeta(x)$, where $\zeta \in W^{1,2}_0(\Omega)\cap L^{\infty}(\Omega)$ and $\psi \in C^{0,1}[0,T]$. In the following, we fix a specific choice of the time-dependent function $\psi$.
For $t\in (0,T)$, $\varepsilon>0$, and $\delta\in (0,\varepsilon)$ such that $t+\varepsilon+\delta<T$, we set
$$
\psi(\tau):= \left\{ 
\begin{aligned}
&0 &&\textrm{for } \tau\in (0,t) \cap (t+\varepsilon+\delta, \infty),\\
&1 &&\textrm{for } \tau\in (t+\delta, t+\varepsilon),\\
&\textrm{linear } &&\textrm{for } \tau\in (t, t+\delta)\textrm{ and } \tau\in (t+\varepsilon, t+\varepsilon+\delta). 
\end{aligned} \right.
$$
With this choice of $\psi$, relation \eqref{eeB} yields 
\begin{equation*}
    \begin{aligned}
        \frac{1}{\delta}\int_{t+\varepsilon}^{t+\varepsilon+\delta} \int_{\Omega}\eta \zeta \dx \dtau - \frac{1}{\delta}\int_{t}^{t+\delta}\int_{\Omega}\eta \zeta\dx\dtau - \int_{t}^{t+\varepsilon+\delta}  
        \intom \eta \bu \cdot \nabla \zeta\, \psi
        \dx \dtau +  \int_{t}^{t+\varepsilon+\delta} \intom \kappa(\vth) \nabla\eta \cdot \nabla \zeta\, \psi \dx \dtau \\
    =  \int_{t}^{t+\varepsilon+\delta} \intom \frac{g\,  \zeta\, \psi }{\vth} \dx \dtau +  \int_{t}^{t+\varepsilon+\delta} \intom \kappa(\vth)
    |\nabla \eta|^2 \zeta \psi \dx \dtau.
    \end{aligned} 
\end{equation*}
Letting $\delta\to 0$ in the above identity, we deduce that for almost all $t\in (0,T)$, all $\varepsilon>0$ such that $t+\varepsilon<T$, and all $\zeta\in W^{1,2}_0(\Omega)\cap L^{\infty}(\Omega)$ it holds
\begin{equation*}
    \begin{aligned}
        \int_{\Omega}(\eta(t+\varepsilon)-\eta(t)) \zeta \dx  - \int_{t}^{t+\varepsilon}  
        \intom \eta \bu \cdot \nabla \zeta
        \dx \dtau +  \int_{t}^{t+\varepsilon} \intom \kappa(\vth) \nabla\eta \cdot \nabla \zeta \dx \dtau \\
    =  \int_{t}^{t+\varepsilon} \intom \frac{g\,  \zeta }{\vth} \dx \dtau +  \int_{t}^{t+\varepsilon} \intom \kappa(\vth)
    |\nabla \eta|^2 \zeta  \dx \dtau,
    \end{aligned} 
\end{equation*}
which can be rewritten as 
\begin{equation}
\label{eeB1}
    \begin{aligned}
        \int_{\Omega}\partial_t \eta_{\varepsilon}(t) \zeta \dx  - \fint_{t}^{t+\varepsilon}  
        \intom \eta \bu \cdot \nabla \zeta
        \dx \dtau +  \fint_{t}^{t+\varepsilon} \intom \kappa(\vth) \nabla\eta \cdot \nabla \zeta \dx \dtau \\
    =  \fint_{t}^{t+\varepsilon} \intom \frac{g\,  \zeta }{\vth} \dx \dtau +  \fint_{t}^{t+\varepsilon} \intom \kappa(\vth)
    |\nabla \eta|^2 \zeta  \dx \dtau,
    \end{aligned} 
\end{equation}
where we defined 
\begin{align}\label{etahlp}
    \eta_\varepsilon(t,x) = \fint_t^{t+\varepsilon} \eta(\tau,x) \dtau.
\end{align} 
Note that $\eta_{\varepsilon} \in C([0,3T/4]; L^q(\Omega))$, for all $\varepsilon\in(0,T/4)$ and $q \in [1,\infty)$. This is a consequence of the definition \eqref{etahlp}, the identity $\eta = \ln \vth$ together with the uniform lower bound $\vth \ge \mu > 0$, and condition \eqref{vb3}. To simplify the notation, we denote $h:=\frac{g}{\vth} + \kappa(\theta) |\nabla \eta|^2$ and note that $h\in L^1(Q)$.  Next, we show that $\eta_{\varepsilon}$ is a Cauchy sequence in $C([0,3T/4]; L^1(\Omega))$ as $\varepsilon \to 0$, and consequently $\eta \in \mathcal{C}([0,3T/4]; L^1(\Omega))$.

To show that $\eta_{\varepsilon}$ is Cauchy, take $\varepsilon_1, \varepsilon_2\in(0,T/8)$ arbitrary. Then it follows from \eqref{eeB1} that for $t\in(0,3T/4)$
\begin{equation}
\label{eeBC}
    \begin{aligned}
        \int_{\Omega}\partial_t (\eta_{\varepsilon_1}(t)-\eta_{\varepsilon_2}(t)) \zeta \dx &= \fint_{t}^{t+\varepsilon_1}  
        \intom \eta \bu \cdot \nabla \zeta \dx \dtau -\fint_{t}^{t+\varepsilon_2} \intom \eta \bu \cdot \nabla \zeta \dx \dtau\\
       &+\fint_{t}^{t+\varepsilon_2} \intom \kappa(\vth) \nabla\eta \cdot \nabla \zeta \dx \dtau - \fint_{t}^{t+\varepsilon_1} \intom \kappa(\vth) \nabla\eta \cdot \nabla \zeta \dx \dtau \\
     &+\fint_{t}^{t+\varepsilon_1} \intom h\, \zeta\dx \dtau-\fint_{t}^{t+\varepsilon_2} \intom h\, \zeta \dx \dtau.
    \end{aligned} 
\end{equation}
We set $\zeta := \mathcal{T}_1(\eta_{\varepsilon_1}(t) - \eta_{\varepsilon_2}(t))$ in \eqref{eeBC}.
Due to the definition of $\mathcal{T}_1$, see \eqref{tk}, and the assumption \eqref{vb4}, we have
$$
\|\zeta\|_{L^{\infty}(Q_7)} + \int_0^{7T/8} \|\zeta\|_{1,2}^2 \le C,
$$
where $Q_7=(0,7T/8)\times\Omega$ and $C$ is independent of $\varepsilon$, so such $\zeta$ is indeed a suitable choice for the test function. 
By the properties of the time mollification \eqref{etahlp}, we can fix a time $t^*\in (3T/4,7T/8)$ such that 
\begin{equation}
\eta_{\varepsilon}(t^*) \to \eta(t^*) \textrm{ strongly in } L^2(\Omega). \label{hlptime}
\end{equation}

Integrating \eqref{eeBC} over $t \in (\ell,t^*)$, with $\ell\in(0,3T/4)$ and using the H\"{o}lder inequality together with the uniform bound above, we deduce that
\begin{equation}
\label{eeBD}
    \begin{aligned}
        &-\int^{t^*}_{\ell}\int_{\Omega}\partial_t (\eta_{\varepsilon_1}(t)-\eta_{\varepsilon_2}(t)) \mathcal{T}_1(\eta_{\varepsilon_1}(t)-\eta_{\varepsilon_2}(t)) \dx \dt \\
        &= -\int^{t^*}_{\ell} \intom \left(\fint_{t}^{t+\varepsilon_1}  
         \eta \bu \dtau -  \fint_{t}^{t+\varepsilon_2}\eta \bu \dtau \right) \cdot \nabla \mathcal{T}_1(\eta_{\varepsilon_1}(t)-\eta_{\varepsilon_2}(t)) \dx \dt\\
       &-\int^{t^*}_{\ell} \intom \left(\fint_{t}^{t+\varepsilon_2} \kappa(\vth) \nabla\eta \dtau -\fint_{t}^{t+\varepsilon_1}  \kappa(\vth) \nabla\eta \dtau \right) \cdot \nabla \mathcal{T}_1(\eta_{\varepsilon_1}(t)-\eta_{\varepsilon_2}(t))  \dx \dt\\
     &-\int^{t^*}_{\ell} \intom\left(  \fint_{t}^{t+\varepsilon_1}  h \dtau -\fint_{t}^{t+\varepsilon_2} h \dtau \right)  \mathcal{T}_1(\eta_{\varepsilon_1}(t)-\eta_{\varepsilon_2}(t))\dx \dt \\
     &\le C\bigg(\left\|\fint_{t}^{t+\varepsilon_1}  
         \eta \bu \dtau -  \fint_{t}^{t+\varepsilon_2}\eta \bu \dtau \right\|_{L^2(Q_7)}+\left\|\fint_{t}^{t+\varepsilon_2} \kappa(\vth) \nabla\eta \dtau -\fint_{t}^{t+\varepsilon_1}  \kappa(\vth) \nabla\eta \dtau \right\|_{L^2(Q_7)}\\
     &\qquad +\left\|  \fint_{t}^{t+\varepsilon_1}  h \dtau -\fint_{t}^{t+\varepsilon_2} h \dtau \right\|_{L^1(Q_7)}\bigg).
    \end{aligned} 
\end{equation}
 From the properties of the Bochner integral, it follows that the quantities 
 $
 \fint_{t}^{t+\varepsilon}  
         \eta \bu \dtau
 $,
 $
 \fint_{t}^{t+\varepsilon} \kappa(\vth) \nabla\eta \dtau
 $,
 and
 $
 \fint_{t}^{t+\varepsilon_1}  h \dtau
 $
 form Cauchy sequences in $L^2(Q_7)$ and $L^1(Q_7)$, respectively.
 Consequently, the right-hand side of \eqref{eeBD} can be arbitrarily small if considering small $\epsilon_1$ and $\epsilon_2$. Recall the definition of $\mathcal{G}_k$ given in \eqref{dfGk}. Then we may rewrite
$$
-\int^{t^*}_{\ell}\!\!\int_{\Omega} \partial_t (\eta_{\varepsilon_1}(t)-\eta_{\varepsilon_2}(t))\,
\mathcal{T}_1(\eta_{\varepsilon_1}(t)-\eta_{\varepsilon_2}(t)) \, dx\, dt
=
\int_{\Omega} \mathcal{G}_1(\eta_{\varepsilon_1}(\ell)-\eta_{\varepsilon_2}(\ell))\, dx
-
\int_{\Omega} \mathcal{G}_1(\eta_{\varepsilon_1}(t^*)-\eta_{\varepsilon_2}(t^*))\, dx.
$$
Due to the choice of $t^*$ in \eqref{hlptime}, the sequence $\eta_{\varepsilon}(t^*)$ is Cauchy in $L^1(\Omega)$. Together with the identity above, estimate \eqref{eeBD} implies that for every $\delta>0$, there exists $\varepsilon_0\in(0,t/4)$ such that, for all $\varepsilon_1,\varepsilon_2 \le \varepsilon_0$, it holds
$$
\sup_{t\in (0,3T/4)} \|\mathcal{G}_1(\eta_{\varepsilon_1}(t)-\eta_{\varepsilon_2}(t))\|_{1} \le \delta.
$$
Moreover, the definition \eqref{dfGk} of $\mathcal{G}_1$, its convexity and Jensen's inequality imply
$$
\sup_{t\in (0,3T/4)} \|\eta_{\varepsilon_1}(t)-\eta_{\varepsilon_2}(t)\|_{1} \le (\mathcal{G}_1)_{-1}(\delta).
$$
Hence $\eta_{\varepsilon}$ is a Cauchy sequence in $C([0,3T/4]; L^1(\Omega))$ which gives a representative of $\eta$ in $C([0,3T/4]; L^1(\Omega))$. Defining $\tilde{\eta}_\varepsilon(t)=\eta_\varepsilon(t-\varepsilon)$ for $t\in(0,T)$ such that $t-\varepsilon\in(0,T)$ 
one can mimic the procedure starting from \eqref{eeB1} to get a representative of $\eta$ in $C([T/4,T]; L^1(\Omega))$. Consequently, there is a representative $\eta\in\mathcal{C}([0,T]; L^1(\Omega))$. It follows from \eqref{eeB} that for this representative $\eta(0)=\eta_0$ in $L^2(\Omega)$. 

Furthermore, since $\sup_{t\in (0,3T/4)} \|\eta_{\varepsilon}(t)\|_q + \sup_{t\in (T/4,T)} \|\tilde\eta_{\varepsilon}(t)\|_q \le C$, with $C>0$ independent of $\epsilon\in(0,T/4)$, for all $q \in [1,\infty)$, a standard Lebesgue interpolation argument yields the desired result
\begin{equation}
\eta(0)=\eta_0,\qquad\eta \in \mathcal{C}([0,T]; L^q(\Omega))
\qquad \text{for all } q\in [1,\infty).
\label{Ceta}
\end{equation}

We now establish the continuity in time of the temperature. Recall the definition \eqref{tk} and, for the sake of brevity, set
$$
    \vth^{m,\delta}:= \mathcal{T}_{m,\delta}(\vth)=\mathcal{T}_{m,\delta}(e^\eta).
$$
Since $\mathcal{T}_{m,\delta}\circ \exp$ is Lipschitz function on $\rr$
we obtain the estimate
$$
|\vth^{m,\delta}(t,x)-\vth^{m,\delta}(\tau,y)|\le C(m,\delta)|\eta(t,x) - \eta(\tau,y)|.
$$
Combining this bound with \eqref{Ceta} yields
\begin{equation}\label{Tkc}
\vth^{m,\delta}(0)=\mathcal{T}_{m,\delta}(\vth_0),\qquad
\vth^{m,\delta} \in C([0,T]; L^q(\Omega)) \qquad \textrm{for all } q\in [1,\infty), \; m>0, \; \delta\in(0,m).
\end{equation}
Similarly, since $\mathcal{T}_{m,\delta}\circ \operatorname{id}^\alpha\circ\exp$ is Lipschitz,
we obtain 
\begin{equation}\label{Tkc2}
\mathcal{T}_{m,\delta} (\vth^{\alpha}) \in C([0,T]; L^q(\Omega)) \qquad \textrm{for all } q\in [1,\infty), \; m>0, \; \delta\in(0,m), \;\alpha \in [0,1].
\end{equation}
Our goal is to determine under which conditions we can remove the truncation function $\mathcal{T}_{m,\delta}$ in the relations \eqref{Tkc} and \eqref{Tkc2}. We first introduce the time mollifications for $t,t+\varepsilon\in(0,T)$ and $\varepsilon\in(0,T/4)$
\begin{equation*}
\begin{aligned}
    \vth^{m,\delta}_\varepsilon(t,x) &= \fint_t^{t+\varepsilon} \mathcal{T}_{m,\delta}(\vth(\tau,x)) \dtau, \qquad && \vth_\varepsilon(t,x) = \fint_t^{t+\varepsilon} \vth(\tau,x) \dtau,\\
    \vth^{m,\delta,\alpha}_\varepsilon(t,x) &= \fint_t^{t+\varepsilon} \mathcal{T}_{m,\delta}(\vth^{\alpha}(\tau,x)) \dtau, \qquad && \vth^{\alpha}_\varepsilon(t,x) = \fint_t^{t+\varepsilon} \vth^{\alpha}(\tau,x) \dtau.
\end{aligned}
\end{equation*}
Since $\vth^{\alpha}_\varepsilon, \, \vth_\varepsilon \in C([0,3T/4];L^1(\Omega))$, if we show that these sequences are Cauchy as $\varepsilon \to \infty$, the continuity of $\vth$ and $\vth^\alpha$ also holds. By \eqref{Tkc}--\eqref{Tkc2}, we know that $\vth^{m,\delta}_\varepsilon$ and $\vth^{m,\delta,\alpha}_\varepsilon$ are Cauchy for any fixed $m,\delta$. We thus use the triangle inequality to obtain for $0<\delta\leq m/2$ the estimate
\begin{equation}\label{TOSO1}
\begin{split}
\|\vth^{\alpha}_{\varepsilon_1}(t)-\vth^{\alpha}_{\varepsilon_2}(t)\|_1 &\le \|\vth^{m,\delta,\alpha}_{\varepsilon_1}(t)-\vth^{m,\delta,\alpha}_{\varepsilon_2}(t)\|_1+\|\vth^{\alpha}_{\varepsilon_1}(t)-\vth^{m,\delta,\alpha}_{\varepsilon_1}(t)\|_1
+\|\vth^{\alpha}_{\varepsilon_2}(t)-\vth^{m,\delta,\alpha}_{\varepsilon_2}(t)\|_1\\
&\le \|\vth^{m,\delta,\alpha}_{\varepsilon_1}(t)-\vth^{m,\delta,\alpha}_{\varepsilon_2}(t)\|_1+\fint_t^{t+\varepsilon_1}\!\!\!\int_{\{\vth^{\alpha}>m/2\}}\!\!\! |\vth^{\alpha}(\tau)|\dx \dtau +\fint_t^{t+\varepsilon_2}\!\!\!\int_{\{\vth^{\alpha}>m/2\}} \!\!\!|\vth^{\alpha}(\tau)|\dx \dtau\\
&\le \|\vth^{m,\delta,\alpha}_{\varepsilon_1}(t)-\vth^{m,\delta,\alpha}_{\varepsilon_2}(t)\|_1+2\essup_{t\in (0,T)}\int_{\{\vth^{\alpha}>m/2\}} |\vth^{\alpha}(t)|\dx \\
&\le \|\vth^{m,\delta,\alpha}_{\varepsilon_1}(t)-\vth^{m,\delta,\alpha}_{\varepsilon_2}(t)\|_1+2\left(\frac m2\right)^{\alpha-1}\essup_{t\in (0,T)}\int_{\Omega} |\vth(t)|\dx
\end{split}
\end{equation}
and
\begin{equation}\label{TOSO2}
\begin{split}
\|\vth_{\varepsilon_1}(t)-\vth_{\varepsilon_2}(t)\|_1 &\le \|\vth^{m,\delta}_{\varepsilon_1}(t)-\vth^{m,\delta}_{\varepsilon_2}(t)\|_1+\|\vth_{\varepsilon_1}(t)-\vth^{m,\delta}_{\varepsilon_1}(t)\|_1
+\|\vth_{\varepsilon_2}(t)-\vth^{m,\delta}_{\varepsilon_2}(t)\|_1\\
&\le \|\vth^{m,\delta}_{\varepsilon_1}(t)-\vth^{m,\delta}_{\varepsilon_2}(t)\|_1+\fint_t^{t+\varepsilon_1}\!\!\!\int_{\{\vth>m/2\}} \!\!\!|\vth(\tau)|\dx \dtau +\fint_t^{t+\varepsilon_2}\!\!\!\int_{\{\vth>m/2\}} \!\!\!|\vth(\tau)|\dx \dtau\\
&\le \|\vth^{m,\delta}_{\varepsilon_1}(t)-\vth^{m,\delta}_{\varepsilon_2}(t)\|_1+2\essup_{t\in (0,T)}\int_{\{\vth>m/2\}} |\vth(t)|\dx.
\end{split}
\end{equation}
To prove \eqref{contir}, we use the estimate \eqref{TOSO1} with $\alpha\in(0,1)$ together with \eqref{vb2B}. We obtain
\begin{equation*}
\begin{split}
\|\vth^{\alpha}_{\varepsilon_1}-\vth^{\alpha}_{\varepsilon_2}\|_{C([0,3T/4]; L^1(\Omega))} 
&\le \|\vth^{m,\delta,\alpha}_{\varepsilon_1}-\vth^{m,\delta,\alpha}_{\varepsilon_2}\|_{C([0,3T/4]; L^1(\Omega))} 
    + Cm^{\alpha-1}.
\end{split}
\end{equation*}
The second term on the right-hand side is independent of $\varepsilon$ and can be made arbitrarily small by a suitable choice of $m$. Moreover, for any fixed $m$, the family $\vth^{m,\delta,\alpha}_\varepsilon$ forms a Cauchy sequence. It follows that $\vth^{\alpha}_\varepsilon$ is also Cauchy and we get $\vth^\alpha \in C([0,3T/4];L^1(\Omega))$. Since we can proceed similarly also on $[T/4,T]$ we obtain \eqref{contir}. 

The above procedure cannot be used in the case $\alpha = 1$. Here, we require additional uniform equiintegrability assumption:
\begin{equation}\label{equi}
\lim_{m\to +\infty} \essup_{t\in (0,T)} \int_{\Omega} \vth(t,x)\chi_{\{|\vth(t,x)|\ge m\}}\dx=0.
\end{equation}
Under this assumption, inequality \eqref{TOSO2} implies continuity $\vth\in C([0,T]; L^1(\Omega))$. We shall later verify that assumption \eqref{equi} can in fact be derived from \eqref{Equi1}.

\bigskip

To prove the remaining results, we begin by showing that equation \eqref{iebR} is satisfied by $\vth^{m,\delta}$. We test \eqref{eeB1} by $e^{\eta_\varepsilon}\,\mathcal{T}'_{m,\delta}(e^{\eta_\varepsilon})\,\varphi$, where $\varphi \in C^1_0((-\infty,T)\times \Omega)$ is arbitrary and $\varepsilon>0$ is sufficiently small so that the resulting terms are well defined in $(0,T)$. Integrating the result over the time interval $(0,T)$ and performing integration by parts, we arrive at
\begin{equation*}
    \begin{aligned}
        &-\int_0^T \int_{\Omega}\mathcal{T}_{m,\delta}(e^{\eta_\varepsilon})\,\partial_t\varphi \dx \dt  +\int_0^T\fint_{t}^{t+\varepsilon}  
        \intom \nabla \eta \cdot \bu \left(e^{\eta_\varepsilon}\,\mathcal{T}'_{m,\delta}(e^{\eta_\varepsilon})\,\varphi\right)
        \dx \dtau \dt\\
        &\qquad +  \int_0^T \fint_{t}^{t+\varepsilon} \intom \kappa(\vth) \nabla\eta \cdot \nabla \left(e^{\eta_\varepsilon}\,\mathcal{T}'_{m,\delta}(e^{\eta_\varepsilon})\,\varphi\right) \dx \dtau \dt\\
    &=  \int_0^T \fint_{t}^{t+\varepsilon} \intom \left(\frac{g}{\vth} + \kappa(\vth)
    |\nabla \eta|^2 \right)e^{\eta_\varepsilon}\,\mathcal{T}'_{m,\delta}(e^{\eta_\varepsilon})\,\varphi \dx \dtau \dt
    +\int_{\Omega}\mathcal{T}_{m,\delta}(e^{\eta_\varepsilon(0)})\varphi(0) \dx.
    \end{aligned} 
\end{equation*}
Due to the assumptions \eqref{vb1} and
\eqref{vb4}, the properties \eqref{propertiestke} of $\mathcal{T}_{m,\delta}'$, \eqref{Tkc}, and substituting $e^{\eta}=\vth$ from the definition of $\eta$, we can pass to the limit $\varepsilon\to 0$ in the above identity to deduce
\begin{equation}
\label{eeB1.1}
    \begin{aligned}
        &-\int_0^T \int_{\Omega}\mathcal{T}_{m,\delta}(\vth)\,\partial_t\varphi \dx \dt  +\int_0^T  
        \intom \nabla\mathcal{T}_{m,\delta}(\vth) \cdot \bu \,\varphi
        \dx \dt +  \int_0^T \intom \kappa(\vth) \nabla\eta \cdot \nabla \left(\vth \,\mathcal{T}'_{m,\delta}(\vth)\,\varphi\right) \dx \dtau \dt\\
    &\qquad=  \int_0^T\intom \left(\frac{g}{\vth} + \kappa(\vth)
    |\nabla \eta|^2 \right)\vth\,\mathcal{T}'_{m,\delta}(\vth)\,\varphi \dx  \dt
    +\int_{\Omega}\mathcal{T}_{m,\delta}(\vth_0)\varphi(0) \dx.
    \end{aligned} 
\end{equation}
Finally, we can rearrange the above identity so that for all $\varphi \in C^1_0((-\infty,T)\times \Omega)$ it holds 
\begin{equation}
\label{eeB1.2}
    \begin{aligned}
        &-\int_0^T \int_{\Omega}\mathcal{T}_{m,\delta}(\vth)\,\partial_t\varphi \dx \dt  +\int_0^T  
        \intom \nabla\mathcal{T}_{m,\delta}(\vth) \cdot \bu \,\varphi
        \dx \dt +  \int_0^T \intom \kappa(\vth) \nabla \vth \cdot  \nabla \varphi \, \mathcal{T}'_{m,\delta}(\vth) \dx  \dt\\
        &\quad =  \int_0^T\intom g\,\mathcal{T}'_{m,\delta}(\vth)\,\varphi \dx  \dt +\int_{\Omega}\mathcal{T}_{m,\delta}(\vth_0) \varphi(0)\dx-    \int_0^T \intom \kappa(\vth) \mathcal{T}''_{m,\delta}(\vth) |\nabla \vth|^2 \varphi \dx \dt,
    \end{aligned} 
\end{equation}
which proves that the partial differential identity \eqref{iebR} is satisfied. 

Let us now show that \eqref{vb2}--\eqref{tls} follow from \eqref{eeB}, the assumptions \eqref{vb2B} and \eqref{vb4}, and the fact that $\hat{\eta} \in W^{1,2}(\Omega)\cap L^{\infty}(\Omega)$. 
In the following, we proceed at a formal level. A fully rigorous argument can be obtained by applying time mollification, in analogy with the treatment of the equation for $\eta$, see \eqref{eeB1}--\eqref{eeBD}. Since this has already been discussed in detail, we do not repeat the argument here.
We set $M := \|\hat{\eta}\|_{L^{\infty}(\partial \Omega)}$ and introduce the test function
$$
\varphi := \frac{(e^{\alpha \eta} - e^{\alpha M})_+}{1+\varepsilon e^{\alpha \eta}},
$$  where $\alpha \in (0,1)$ and $\varepsilon>0$. Observe that $\varphi \in L^{\infty}(Q)\cap L^2(0,T; W^{1,2}_0(\Omega))$. Next, we define
$$
G_{\varepsilon}(s) := \int_0^s \frac{(e^{\alpha t} - e^{\alpha M})_+}{1+\varepsilon e^{\alpha t}} \,\dt
$$
and note the elementary bounds
\begin{align}\label{boundstest}
    |\varphi| \le C\vth && |G_{\varepsilon}(\eta)| \le C\vth.
\end{align}
We now use $\varphi$ as a test function in \eqref{eeB}. Integrating by parts, exploiting the identity $\diver \bu = 0$, and using \eqref{Ceta}, 
we obtain
\begin{equation*}
    \begin{aligned}
        &\int_{\Omega} G_{\varepsilon}(\eta(T))-G_{\varepsilon}(\eta(0)) \dx  + \intt \intom \frac{\kappa(\vth) |\nabla\eta|^2 \alpha e^{\alpha \eta}}{{1+\varepsilon e^{\alpha \eta}}} \chi_{\{\eta>M\}} \dx \dt - \intt \intom \frac{g (e^{\alpha \eta} - e^{\alpha M})_+}{\vartheta(1+\varepsilon e^{\alpha \eta})} \dx \dt\\
    & = \intt \intom \kappa(\vth)
    |\nabla \eta|^2 \frac{(e^{\alpha \eta} - e^{\alpha M})_+}{1+\varepsilon e^{\alpha \eta}} +\kappa(\vth)
    |\nabla \eta|^2 \frac{(e^{\alpha \eta} - e^{\alpha M})_+ \varepsilon \alpha e^{\alpha \eta}}{(1+\varepsilon e^{\alpha \eta})^2} \dx \dt.
    \end{aligned} 
\end{equation*}
After rearranging the terms, neglecting the nonnegative contributions, and using the bounds \eqref{boundstest}, we obtain the estimate
\begin{equation*}
    \begin{aligned}
     (1-\alpha) \intt \intom \kappa(\vth)
    |\nabla \eta|^2 \frac{e^{\alpha \eta}}{1+\varepsilon e^{\alpha \eta}} \dx \dt
    \le C \int_{\Omega} \vth(T) \dx
    + C \intt \intom \kappa(\vth) |\nabla \eta|^2 e^{\alpha M}  + |g|\dx \dt.
    \end{aligned}
\end{equation*}
Letting $\varepsilon \to 0$, using the fact that $\vth=\ln \eta$ and the assumption \eqref{kap}, we deduce 
\begin{equation}
\label{eeBa}
    \begin{aligned}
     (1-\alpha) \intt \intom \frac{|\nabla \vth|^2}{\vth^{2-\alpha}}\dx \dt
    \le C \int_{\Omega} \vartheta(T) \dx
    + C \intt \intom  |\nabla \eta|^2 e^{\alpha M}  + |g|\dx \dt \le C.
    \end{aligned}
\end{equation}
This estimate yields \eqref{vb2}. The bounds \eqref{vb3} and \eqref{tls} then follow by interpolation, combined with assumption \eqref{vb2B}, see \cite{BoGa89,Pr97,BuCoMa11}.
 
%

Moreover, if $p > (d+2)/2$, it follows from \eqref{vb1} and \eqref{vb3} that $\bu \vth \in L^1(Q;\mathbb{R}^d)$. In combination with the bound \eqref{tls} for the temperature gradient, this allows us to pass to the limit $m \to \infty$ in \eqref{eeB1.2}. An integration by parts in the second integral then yields
\begin{equation}
\label{eeB1.3}
    \begin{aligned}
        &-\int_0^T \int_{\Omega}\vth\,\partial_t\varphi \dx \dt  -\int_0^T  
        \intom \vth  \bu \cdot \nabla \varphi
        \dx \dt +  \int_0^T \intom \kappa(\vth) \nabla \vth \cdot  \nabla \varphi  \dx  \dt\\
        &\quad=  \int_0^T\intom g\,\varphi \dx  \dt +\int_{\Omega}\vth_0 \varphi(0)\dx -\lim_{m\to \infty} \int_0^T \intom \kappa(\vth) \mathcal{T}''_{m,\delta}(\vth) |\nabla \vth|^2 \varphi \dx \dt.
    \end{aligned} 
\end{equation}
Since the last limit is nonpositive for nonnegative $\varphi$ by \eqref{propertiestke}, we infer that \eqref{iebB} holds. 


\bigskip

It remains to prove the equivalence of statements {\bf A)} -- {\bf E)} and their relation to continuity. 
We begin by proving that \eqref{Equi1} implies that \eqref{Equi2} holds with equality, i.e., that {\bf A)} implies {\bf C)}. Since a rigorous procedure for dealing with time derivatives has already been provided when deriving the equation \eqref{iebR} for $\mathcal{T}_{m,\delta}(\vth)$ from \eqref{eeB} we shall proceed here at a formal level\footnote{One could mollify with respect to the time variable and test by the corresponding quantity.}.
We choose $m$ so that $2\norm{\htet}{\infty}{ } < m$, $k=1$ and $\delta \in (0, m/2)$, and set 
$$
\varphi := \mathcal{T}_{k}(\mathcal{T}_{m,\delta}(\vth) - \psi) \chi_{\{t\in [0,\tau]\}}
$$ 
in \eqref{iebR}, where $\psi \in L^{\infty}(Q) \cap L^2(0,T; W^{1,2}(\Omega)) \cap W^{1,1}(0,T; L^1(\Omega))$ and satisfies $\psi = \hat{\vth}$ on $\partial \Omega$. We obtain
\begin{equation}
\label{eeB1.220}
    \begin{aligned}
        &\int_0^{\tau} \int_{\Omega}\partial_t \mathcal{T}_{m,\delta}(\vth)\mathcal{T}_{k}(\mathcal{T}_{m,\delta}(\vth) - \psi) \dx \dt  +\int_0^{\tau}  
        \intom \nabla\mathcal{T}_{m,\delta}(\vth) \cdot \bu \mathcal{T}_{k}(\mathcal{T}_{m,\delta}(\vth) - \psi)
        \dx \dt \\
        &\qquad +  \int_0^{\tau} \intom \kappa(\vth)   \nabla  \mathcal{T}_{m,\delta}(\vth)\cdot \nabla \mathcal{T}_{k}(\mathcal{T}_{m,\delta}(\vth) - \psi) \dx  \dt\\
        &=  \int_0^{\tau}\intom g\,\mathcal{T}'_{m,\delta}(\vth)\,\mathcal{T}_{k}(\mathcal{T}_{m,\delta}(\vth) - \psi)\dx  \dt x-    \int_0^{\tau} \intom \kappa(\vth) \mathcal{T}''_{m,\delta}(\vth) |\nabla \vth|^2 \mathcal{T}_{k}(\mathcal{T}_{m,\delta}(\vth) - \psi) \dx \dt.
    \end{aligned} 
\end{equation}
Using integration by parts and the identity $\diver \bu =0$, the first two terms in \eqref{eeB1.220} can be for a.e. $\tau\in(0,T)$ rewritten as
$$
\begin{aligned}
&\int_0^{\tau}\int_{\Omega}\partial_t \mathcal{T}_{m,\delta}(\vth)\mathcal{T}_{k}(\mathcal{T}_{m,\delta}(\vth) - \psi) \dx \dt  +\int_0^{\tau}  
        \intom \nabla\mathcal{T}_{m,\delta}(\vth) \cdot \bu \mathcal{T}_{k}(\mathcal{T}_{m,\delta}(\vth) - \psi)
        \dx \dt\\
&=\int_0^{\tau} \int_{\Omega}\partial_t (\mathcal{T}_{m,\delta}(\vth)-\psi)\mathcal{T}_{k}(\mathcal{T}_{m,\delta}(\vth) - \psi) \dx \dt+\int_0^{\tau} \int_{\Omega}\partial_t\psi \mathcal{T}_{k}(\mathcal{T}_{m,\delta}(\vth) - \psi) \dx \dt  \\
&\qquad +\int_0^{\tau}  
        \intom \nabla(\mathcal{T}_{m,\delta}(\vth)-\psi) \cdot \bu \mathcal{T}_{k}(\mathcal{T}_{m,\delta}(\vth) - \psi)+\int_0^{\tau}  
        \intom \nabla \psi \cdot \bu \mathcal{T}_{k}(\mathcal{T}_{m,\delta}(\vth) - \psi)
        \dx \dt       \\
&= \int_{\Omega}\mathcal{G}_{k}(\mathcal{T}_{m,\delta}(\vth(\tau)) - \psi(\tau)) -\mathcal{G}_{k}(\mathcal{T}_{m,\delta}(\vth_0) - \psi(0)) \dx +\int_0^{\tau} \int_{\Omega}\partial_t\psi \mathcal{T}_{k}(\mathcal{T}_{m,\delta}(\vth) - \psi) \dx \dt  \\
&\qquad +\int_0^{\tau}  
        \intom \nabla \psi \cdot \bu \mathcal{T}_{k}(\mathcal{T}_{m,\delta}(\vth) - \psi)
        \dx \dt.       
\end{aligned}
$$
Furthermore, denoting $M:=\|\psi\|_{L^{\infty}(Q)}$, the third integral in \eqref{eeB1.220} can be rewritten as
$$
\begin{aligned}
 &\int_0^{\tau} \intom \kappa(\vth)   \nabla  \mathcal{T}_{m,\delta}(\vth)\cdot \nabla \mathcal{T}_{k}(\mathcal{T}_{m,\delta}(\vth) - \psi) \dx  \dt\\
 &\quad =\int_0^{\tau} \intom \kappa(\vth)   \nabla  \mathcal{T}_{k+M}(\mathcal{T}_{m,\delta}(\vth))\cdot \nabla \mathcal{T}_{k}(\mathcal{T}_{m,\delta}(\vth) - \psi) \dx  \dt.
\end{aligned}
$$
Substituting into \eqref{eeB1.220}, recalling $\delta\leq m$, 
and letting $m\to \infty$, we arrive at
\begin{equation}
\label{eeB1.2234}
    \begin{aligned}
        &\int_{\Omega}\mathcal{G}_{k}(\vth(\tau) - \psi(\tau)) -\mathcal{G}_{k}(\vth_0 - \psi(0)) \dx +\int_0^{\tau} \int_{\Omega}\partial_t\psi \mathcal{T}_{k}(\vth - \psi) \dx \dt  \\
&\qquad +\int_0^{\tau}  
        \intom \nabla \psi \cdot \bu \mathcal{T}_{k}(\vth - \psi)
        \dx \dt       +  \int_0^{\tau} \intom \kappa(\vth)   \nabla  \mathcal{T}_{k+M}(\vth)\cdot \nabla \mathcal{T}_{k}(\vth - \psi) \dx  \dt\\
        &=  \int_0^{\tau}\intom g\,\mathcal{T}_{k}(\vth - \psi)\dx  \dt x-   \lim_{m\to \infty} \int_0^{\tau} \intom \kappa(\vth) \mathcal{T}''_{m,m/2}(\vth) |\nabla \vth|^2 \mathcal{T}_{k}(\mathcal{T}_{m,m/2}(\vth) - \psi) \dx \dt.
    \end{aligned} 
\end{equation}
Using \eqref{Equi1} to identify the limit, we obtain for a.e. $\tau\in(0,T)$ \eqref{Equi2} with equality, which is  {\bf C)}. 

\bigskip

Since {\bf C)} clearly implies {\bf B)}, we want to show that {\bf B)} $\implies$ {\bf A)} to have the equivalence {\bf A)}--{\bf C)}. Suppose that there exists a solution $\vth_1$ satisfying \eqref{Equi2}. We then define $\kappa_1(t,x) := \kappa(\vth_1(t,x))$ and construct a solution to an auxiliary problem with prescribed space- and time-dependent heat conductivity $\kappa_1$. We first show that the constructed solution satisfies \eqref{Equi1}. We then prove that it coincides with $\vth_1$, thereby establishing {\bf A)}.

Hence, let $\vth_1$ be a solution fulfilling \eqref{Equi2} and define $\kappa_1$ as above. We introduce the approximate problem:
\begin{equation}
\begin{split}
\partial_t \vth^n + \bu^n \cdot \nabla \vth^n -\diver (\kappa_1\nabla \vth^n) &= \mathcal{T}_n (g) \quad \textrm{ in } Q,\\
\vth^n&= \hat{\vth} \quad \textrm{ on }(0,T)\times \partial \Omega,\\
\vth^n(0)&=\mathcal{T}_n (\vth_0) \quad \textrm{ in } \Omega,
\end{split}
\label{app2}
\end{equation}
where $\bu^n \in L^{\infty}(Q)\cap L^2(0,T; L^2_{0,\diver})$ fulfilling 
$$
\bu^n \to \bu \textrm{ strongly in }L^2 (Q).
$$
Note that the standard theory gives the existence of $\vth^n\in L^2(0,T;W^{1,2}(\Omega))\cap C([0,T];L^2(\Omega))$ with $\partial_t\vth\in L^2(0,T;(W^{1,2}_0(\Omega))^*$. As in 
Section~\ref{sec:ntn}, we infer that the sequence $\vth^n$ satisfies the uniform bounds in the spaces specified in \eqref{vb2}--\eqref{contir}. Moreover, there exists $\vth \in L^{\infty}(0,T; L^1(\Omega))$ such that, up to a non-relabeled subsequence, for all $k>0$, 
\begin{align}
\vth^n&\to \vth &&\textrm{ a.e. in $Q$},
\label{nt1}\\
\mathcal{T}_k(\vth^n) &\rightharpoonup \mathcal{T}_k(\vth) &&\textrm{ in } L^{2}(0,T; W^{1,2}\om).\label{nt2}
\end{align}
At this moment, it is not known whether $\vth$ is nonnegative, as the nonnegativity of $g$ is not assumed.

First, we show that \(\vartheta\) satisfies \eqref{Equi1}. To this end, we consider an arbitrarily $\tau\in(0,T)$ such that $\vth^n(\tau)\to\vth(\tau)$ in $L^1(\Omega)$. Note that this holds for a.e. $\tau\in(0,T)$ and so $\tau$ can be chosen arbitrarily close to $T$. We test the weak formulation of \eqref{app2} with \((1 - \mathcal{T}'_{m,\delta}(\vartheta^n_+))\chi_{(0,\tau)}\) which belongs to \(L^2(0,T;W^{1,2}_0(\Omega))\) provided \(m \ge 2\|\hat{\vartheta}\|_{\infty}\) and $\delta \in (0, m/2)$. Consequently, we formally obtain 
$$
\begin{aligned}
0&=\int_0^\tau \intom \partial_t \vth^n (1 - \mathcal{T}'_{m,\delta}(\vth^n_+)) + \bu^n \cdot \nabla \vth^n (1 - \mathcal{T}'_{m,\delta}(\vth^n_+)) -\kappa_1|\nabla \vth^n_+|^2 \mathcal{T}''_{m,\delta}(\vth^n_+)- \mathcal{T}_n (g)(1 - \mathcal{T}'_{m,\delta}(\vth^n_+)) \dx \dt\\
&= \intom   \vth^n_+(\tau)-\mathcal{T}_{m,\delta}(\vth^n_+(\tau))-(\mathcal{T}_n(\vth_0) - \mathcal{T}_{m,\delta}(\mathcal{T}_n(\vth_0)))\dx  - \int_0^\tau \intom \kappa_1|\nabla \vth^n_+|^2 \mathcal{T}''_{m,\delta}(\vth^n_+)- |g|(1 - \mathcal{T}'_{m,\delta}(\vth^n_+)) \dx \dt.
\end{aligned}
$$
We now let $n \to \infty$.  We first realize that it follows from the previous equality that the limit
$$
\lim_{n\to+\infty}\int_0^\tau \intom \kappa_1|\nabla \vth^n_+|^2 \mathcal{T}''_{m,\delta}(\vth^n_+)\dx \dt
$$
exists by the choice of $\tau$ and \eqref{nt1}--\eqref{nt2}. By the weak lower semicontinuity of the $L^2$-norm, the pointwise convergence $\vth^n \to \vth$ in $Q$, and the nonpositivity of $\mathcal{T}''_{m,\delta}$ on $(0,+\infty)$, we infer
\begin{equation}
\label{chceme}
\begin{aligned}
-\int_0^\tau \intom |\nabla \vth_+|^2 \mathcal{T}''_{m,\delta}(\vth_+) \dx \dt&\le -\lim_{n\to \infty}\int_0^\tau \intom |\nabla \vth^n_+|^2 \mathcal{T}''_{m,\delta}(\vth^n_+) \dx \dt \\
&\le  C\intom   (\vth_0 - m/2)_+\dx  + C\int_0^T |g| \chi_{\{\vth\ge m/2\}}\dx \dt.
\end{aligned}
\end{equation}
Using Fatou lemma, we can replace $\tau$ with $T$ in the integral on the left-hand side of \eqref{chceme}. Observe that the right-hand side 
vanishes as $m \to \infty$.  Using the properties \eqref{propertiestke} of $\mathcal{T}_{m,\delta}$, we conclude that \eqref{Equi1} is valid for $\vth$.

Repeating the estimate leading to \eqref{chceme} with the test function \((1 - \mathcal{T}'_{m,\delta}(\vartheta^n_-))\chi_{(0,\tau)}\) we get 
\begin{equation}
\label{chceme1}
\begin{aligned}
\lim_{n\to \infty}\int_0^\tau \intom |\nabla \vth^n_-|^2 \mathcal{T}''_{m,\delta}(\vth^n_-) \dx \dt 
\le  C\int_0^T |g| \chi_{\{\vth\le -m/2\}}\dx \dt,
\end{aligned}
\end{equation}
which together with \eqref{chceme} gives 
\begin{equation}
\label{chceme3}
\begin{aligned}
\lim_{n\to \infty}\int_0^\tau \intom |\nabla \vth^n|^2 |\mathcal{T}''_{m,\delta}(\vth^n)| \dx \dt \le  C\intom   (\vth_0 - m/2)_+\dx  + C\int_0^T |g| \chi_{\{|\vth|\ge m/2\}}\dx \dt.
\end{aligned}
\end{equation}

Finally, we show that $\vth=\vth_1$ in $Q$, which implies that $\vth_1$ also satisfies \eqref{Equi1} and $\vth$ is actually nonnegative. 
We mollify $\vth^n$ as 
$$
\vth^{n,m,\delta}_\varepsilon(t,x) = \fint_t^{t+\varepsilon} \mathcal{T}_{m,\delta}(\vth^n(\tau,x)) \, \dtau.
$$
We then test \eqref{Equi2} for $\vth_1$ with $\varphi:= \vth^{n,m,\delta}_\varepsilon$ with $\norm{\htet}{\infty}{ }+1 < m$ and $\delta\in(0,1)$, to obtain for a.e. $\tau\in(0,T-\varepsilon)$
\begin{equation}
\label{Equi3}
\begin{aligned}
& \intom \mathcal{G}_k (\vth_1(\tau) -\vth^{n,m,\delta}_\varepsilon(\tau)) \dx +\int_{0}^{\tau}\intom \partial_t \vth^{n,m,\delta}_\varepsilon\, \mathcal{T}_k(\vth_1-\vth^{n,m,\delta}_\varepsilon) \dx \dt\\
&\qquad +\int_{0}^{\tau} \int_\Omega \nabla \vth^{n,m,\delta}_\varepsilon \cdot  \bu \mathcal{T}_k(\vth_1 - \vth^{n,m,\delta}_\varepsilon)
        \dx \dt + \int_{0}^{\tau} \intom \kappa_1 \nabla\vth_1 \cdot \nabla \mathcal{T}_{k}(\vth_1 - \vth^{n,m,\delta}_\varepsilon) \dx \dt \\
&\leq \int_{0}^{\tau}\intom g \mathcal{T}_{k}(\vth_1 - \vth^{n,m,\delta}_\varepsilon) \dx \dt + \intom \mathcal{G}_k (\vth_0 -\vth^{n,m,\delta}_\varepsilon(0)) \dx.
    \end{aligned} 
\end{equation}  
We now identify the term involving the time derivative of $\vth^{n,m,\delta}_\varepsilon$ in \eqref{Equi3} with the corresponding expression obtained from the approximate problem \eqref{app2}. To avoid unnecessary repetition, we present the argument only at a formal level, since an analogous procedure has already been justified rigorously. More precisely, we multiply \eqref{app2} by $\mathcal{T}_{m,\delta}'(\vth^n)$ and integrate over the time interval $(t, t+\varepsilon)$ with $t\in(0,T-\varepsilon)$. This yields
\begin{equation*}
\begin{split}
\partial_t \vth^{n,m,\delta}_\varepsilon= \fint_t^{t+\varepsilon}-\bu^n \cdot \nabla \mathcal{T}_{m,\delta}(\vth^n)  +\diver (\kappa_1\nabla \mathcal{T}_{m,\delta}(\vth^n)) -\kappa_1 \mathcal{T}_{m,\delta}''(\vth^n)|\nabla \vth^n|^2 + \mathcal{T}_n (g) \mathcal{T}'_{m,\delta}(\vth^n)\ds.
\end{split}
\end{equation*}
Multiplying the above identity by $\mathcal{T}_{k}(\vth_1-\vth^{n,m,\delta}_\varepsilon)$, integrating over $(0,\tau)\times \Omega$, and recalling that $m>\norm{\hat{\vth}}{\infty}{ }+1$ and $\delta\in(0,1)$, we obtain
\begin{equation}
\begin{split}
&\int_0^{\tau}\intom \partial_t \vth^{n,m,\delta}_\varepsilon\mathcal{T}_{k}(\vth_1-\vth^{n,m,\delta}_\varepsilon) \dx \dt \\
&= \int_0^{\tau}\intom \fint_t^{t+\varepsilon}-\bu^n \cdot \nabla \mathcal{T}_{m,\delta}(\vth^n)\ds\mathcal{T}_{k}(\vth_1-\vth^{n,m,\delta}_\varepsilon)  -\fint_t^{t+\varepsilon}\kappa_1\nabla \mathcal{T}_{m,\delta}(\vth^n)\ds\cdot \nabla \mathcal{T}_{k}(\vth_1-\vth^{n,m,\delta}_\varepsilon)\\
&-\fint_t^{t+\varepsilon}\kappa_1 \mathcal{T}_{m,\delta}''(\vth^n)|\nabla \vth^n|^2\ds\mathcal{T}_{k}(\vth_1-\vth^{n,m,\delta}_\varepsilon) + \fint_t^{t+\varepsilon}\mathcal{T}_n (g) \mathcal{T}'_{m,\delta}(\vth^n)\ds \mathcal{T}_{k}(\vth_1-\vth^{n,m,\delta}_\varepsilon)  \dx \dt.
\end{split}
\label{app2rq}
\end{equation}
At this point, we let $\varepsilon \to 0_+$ in \eqref{Equi3}. The integral involving the time derivative is handled by means of identity \eqref{app2rq}, while all remaining terms can be treated directly. Owing to the regularity of the truncation operators and the averaging procedure in time, we may pass to the limit term by term, which yields for a.e. $\tau\in(0,T)$
\begin{equation}
\label{Equi4}
\begin{aligned}
& \intom \mathcal{G}_k (\vth_1(\tau) -\mathcal{T}_{m,\delta}(\vth^{n}(\tau))) \dx  + \int_0^{\tau}\intom   \kappa_1|\nabla\mathcal{T}_{k}(\vth_1-\mathcal{T}_{m,\delta}(\vth^n))|^2\dx \dt \\
&\leq  \int_0^{\tau}\intom \kappa_1 \mathcal{T}_{m,\delta}''(\vth^n)|\nabla \vth^n|^2 \mathcal{T}_{k}(\vth_1-\mathcal{T}_{m,\delta}(\vth^n))  \dx \dt\\
&\qquad +\int_{0}^{\tau} \int_\Omega \nabla \mathcal{T}_{m,\delta}(\vth^n) \cdot  (\bu^n-\bu) \mathcal{T}_k(\vth_1 - \mathcal{T}_{m,\delta}(\vth^n))
        \dx \dt  \\
& \qquad +\int_{0}^{\tau}\intom (g-\mathcal{T}_n (g) \mathcal{T}'_{m,\delta}(\vth^n)) \mathcal{T}_{k}(\vth_1 - \mathcal{T}_{m,\delta}(\vth^n)) \dx \dt + \intom \mathcal{G}_k (\vth_0 -\mathcal{T}_{m,\delta}(\mathcal{T}_n(\vth_0)) \dx\\
&\le  C k\int_0^{\tau}\intom |\mathcal{T}_{m,\delta}''(\vth^n)||\nabla \vth^n|^2 +  |\nabla \mathcal{T}_{m,\delta}(\vth^n)||\bu^n-\bu|+ |g-\mathcal{T}_n (g) \mathcal{T}'_{m,\delta}(\vth^n)| \dx \dt\\
&\qquad  + \intom \mathcal{G}_k (\vth_0 -\mathcal{T}_{m,\delta}(\mathcal{T}_n(\vth_0)) \dx.
    \end{aligned} 
\end{equation}  
We next pass to the limit $n\to \infty$. Realizing that $\tau$ can be chosen such that $\vth^n(\tau)\to\vth(\tau)$ in $L^1(\Omega)$ we can apply \eqref{chceme3} to the previous inequality to  infer that for almost every $\tau\in (0,T)$ the limit function $\vth$ satisfies
\begin{equation*}
\begin{aligned}
& \intom \mathcal{G}_k (\vth_1(\tau) -\mathcal{T}_{m,\delta}(\vth(\tau))) \dx \\
&\le  C k\lim_{n\to \infty} \int_0^{T}\intom |\mathcal{T}_{m,\delta}''(\vth^n)||\nabla \vth^n|^2\dx \dt  + Ck \int_0^T \intom |g-g \mathcal{T}'_{m,\delta}(\vth)| \dx \dt  + \intom \mathcal{G}_k (\vth_0 -\mathcal{T}_{m,\delta}(\vth_0)) \dx\\
&\le  C k\int_0^T \intom |g|\chi_{\{\vth>m\}}  \dx \dt  + \intom \mathcal{G}_k ((\vth_0 -m)_+) + Ck (\vth_0-m)_+\dx.
    \end{aligned} 
\end{equation*}  
Letting $m\to \infty$, we obtain
\begin{equation*}
\begin{aligned}
& \intom \mathcal{G}_k (\vth_1(\tau) -\vth(\tau)) \dx\le 0
\end{aligned} 
\end{equation*}
which yields $\vth=\vth_1$ almost everywhere in $(0,T)\times \Omega$. Since \eqref{Equi1} holds for $\vth$, it also holds for $\vth_1$. Hence, we have established that \eqref{Equi2} implies \eqref{Equi1}. 
Together with the previously shown implications we conclude that all statements {\bf A)} – {\bf C)} are equivalent. 

Now we show that condition \eqref{Equi1} in {\bf A)} implies \eqref{equi}, which in turn yields the continuity property $\vth\in C([0,T]; L^1(\Omega))$ as explained above \eqref{equi}. To this end, let $m > n$ be arbitrary such that $n > \max(\|\hat{\vth}\|_{\infty}+1, 2\|\hat{\vth}\|_{\infty})$. We take equation \eqref{iebR} for $m$ with $\delta_m=m/2$ and for $n$ with $\delta_n\in(0,1)$, and subtract the latter from the former. This yields
\begin{equation}
\label{eeB1.21}
    \begin{aligned}
        &-\int_0^T \int_{\Omega}(\mathcal{T}_{m,\delta_m}(\vth)-\mathcal{T}_{n,\delta_n}(\vth))\,\partial_t\varphi \dx \dt  +\int_0^T  
        \intom \nabla(\mathcal{T}_{m,\delta_m}(\vth)-\mathcal{T}_{n,\delta_n}(\vth)) \cdot \bu \,\varphi
        \dx \dt \\
        &\qquad +    \int_0^T \intom \kappa(\vth) (\mathcal{T}''_{m,\delta_m}(\vth)-\mathcal{T}''_{n,\delta_n}(\vth)) |\nabla \vth|^2 \varphi \dx \dt +  \int_0^T \intom \kappa(\vth) \nabla \vth \cdot  \nabla \varphi \, (\mathcal{T}'_{m,\delta_m}(\vth)-\mathcal{T}'_{n,\delta_n}(\vth)) \dx  \dt\\
        &=  \int_0^T\intom g\,(\mathcal{T}'_{m,\delta_m}(\vth)-\mathcal{T}'_{n,\delta_n}(\vth))\,\varphi \dx  \dt +\int_{\Omega}(\mathcal{T}_{m,\delta_m}(\vth_0) -\mathcal{T}_{n,\delta_n}(\vth_0)) \varphi(0)\dx .
    \end{aligned} 
\end{equation}
Next, since $\mathcal{T}_{m,\delta}(\vth) \in \mathcal{C}([0,T];L^1(\Omega))$, see \eqref{Tkc}, we can deduce from the above identity that also for all $\tau\in [0,T]$ we have
\begin{equation}
\label{eeB1.22}
    \begin{aligned}
        &\int_{\Omega}(\mathcal{T}_{m,\delta_m}(\vth (\tau))-\mathcal{T}_{n,\delta_n}(\vth(\tau)))\,\psi \dx   +\int_0^\tau  
        \intom \nabla (\mathcal{T}_{m,\delta_n}(\vth)-\mathcal{T}_{n,\delta_n}(\vth)) \cdot \bu \,\psi
        \dx \dt \\
        &\qquad +    \int_0^\tau \intom \kappa(\vth) (\mathcal{T}''_{m,\delta_m}(\vth)-\mathcal{T}''_{n,\delta_n}(\vth)) |\nabla \vth|^2 \psi \dx \dt \\
        &\qquad+  \int_0^\tau \intom \kappa(\vth) \nabla \vth \cdot  \nabla \psi \, (\mathcal{T}'_{m,\delta_m}(\vth)-\mathcal{T}'_{n,\delta_n}(\vth)) \dx  \dt\\
        &=  \int_0^\tau\intom g\,(\mathcal{T}'_{m,\delta_m}(\vth)-\mathcal{T}'_{n,\delta_n}(\vth))\psi\dx  \dt +\int_{\Omega}(\mathcal{T}_{m,\delta_m}(\vth_0) -\mathcal{T}_{n,\delta_n}(\vth_0)) \psi\dx 
    \end{aligned} 
\end{equation}
valid for all $\psi \in W^{1,2}_0(\Omega) \cap L^{\infty}(\Omega)$.
We define a sequence of functions \(\psi^k\) by
\[
\psi^k(x) =
\begin{cases}
1, & \text{if } B_{1/k}(x) \subset \Omega,\\[1mm]
0, & \text{if } B_{1/(2k)}(x) \cap \Omega \neq B_{1/(2k)}(x).
\end{cases}
\]
Due to the Lipschitz regularity of \(\Omega\), this can be arranged so that \(|\nabla \psi^k| \le C k\). With this construction, $\psi^k$ belongs to $W^{1,2}_0(\Omega) \cap L^{\infty}(\Omega)$ for all $k$, and \(\psi^k \nearrow 1\) as \(k \to \infty\). Using \(\psi^k\) as a test function in \eqref{eeB1.22}, and noting that
$$
(\mathcal{T}'_{m,\delta_m}(\vth)-\mathcal{T}'_{n,\delta_n}(\vth)) \in L^2(0,T; W^{1,2}_0(\Omega)),
$$
due to the choice of $m,n,\delta_m, \delta_n$,
we deduce by the Hardy inequality
$$
\begin{aligned}
&\left|\int_0^\tau \intom \kappa(\vth) \nabla \vth \cdot  \nabla \psi^k \, (\mathcal{T}'_{m,\delta_m}(\vth)-\mathcal{T}'_{n,\delta_n}(\vth)) \dx  \dt \right| \\
&\quad \le C \int_0^T \|\nabla \mathcal{T}_{2m}(\vth) \chi_{\{|\nabla \psi^k | \neq 0\}}\|_2 \left(\int_{\Omega} |\nabla \psi^k|^2 |\mathcal{T}'_{m,\delta_m}(\vth)-\mathcal{T}'_{n,\delta_n}(\vth)|^2 \dx \right)^{\frac12}  \dt\\
&\quad \le C \int_0^T \|\nabla \mathcal{T}_{2m}(\vth) \chi_{\{|\nabla \psi^k | \neq 0\}}\|_2 \left(\int_{\Omega} \frac{ |\mathcal{T}'_{m,\delta_m}(\vth)-\mathcal{T}'_{n,\delta_n}(\vth)|^2}{(\textrm{dist }(x,\partial \Omega))^2} \dx \right)^{\frac12}  \dt \\
&\quad \le C \int_0^T \|\nabla \mathcal{T}_{2m}(\vth) \chi_{\{|\nabla \psi^k | \neq 0\}}\|_2 \left( \|\nabla \mathcal{T}'_{m,\delta_m}(\vth)\|_2 + \|\nabla\mathcal{T}'_{n,\delta_n}(\vth)\|_2\right) \dt \overset{k \to \infty} \to 0.
\end{aligned}
$$
Hence, by taking $\psi^k$ in \eqref{eeB1.22}, passing to  the limit $k\to\infty$, and applying the above estimate together with integration by parts and the identity $\diver \bu=0$, we deduce that
\begin{equation}
\label{eeB1.23}
    \begin{aligned}
        &\int_{\Omega}\mathcal{T}_{m,\delta_m}(\vth (\tau))-\mathcal{T}_{n,\delta_n}(\vth(\tau)) \dx  +    \int_0^\tau \intom \kappa(\vth) (\mathcal{T}''_{m,\delta_m}(\vth)-\mathcal{T}''_{n,\delta_n}(\vth)) |\nabla \vth|^2  \dx \dt \\
        &=  \int_0^\tau\intom g\,(\mathcal{T}'_{m,\delta_m}(\vth)-\mathcal{T}'_{n,\delta_n}(\vth))\dx  \dt +\int_{\Omega}\mathcal{T}_{m,\delta_m}(\vth_0) -\mathcal{T}_{n,\delta_n}(\vth_0)\dx,
    \end{aligned} 
\end{equation}
is valid for all $\tau \in [0,T]$. 
The next step is to 
let $m\to \infty$. We recall $\delta_m=m/2$.
The term involving \(\mathcal{T}''_{m,m/2}(\vartheta)\) vanishes by the virtue of assumption \eqref{Equi1}. Hence, we deduce that, for almost all \(\tau \in (0,T)\),
\begin{equation*}
    \begin{aligned}
        &\int_{\Omega}\vth (\tau)-\mathcal{T}_{n,\delta_n}(\vth(\tau))\dx   -    \int_0^\tau \intom \kappa(\vth)\mathcal{T}''_{n,\delta_n}(\vth) |\nabla \vth|^2  \dx \dt \\
        &=  \int_0^\tau\intom g\,(1-\mathcal{T}'_{n,\delta_n}(\vth))\dx  \dt +\int_{\Omega}\vth_0 -\mathcal{T}_{n,\delta_n}(\vth_0) \dx.
    \end{aligned} 
\end{equation*}
The statement holds for almost all times \(\tau\),  the  continuity of $\vth$ in $L^1(\Omega)$ with respect to the time variable has not yet been established.
Since the second term on the left-hand side is nonnegative, we may let $\delta_n \to 0_+$ which implies that, for almost all $\tau \in (0,T)$ and for every $n$, it holds
\begin{equation*}
    \begin{aligned}
        &\int_{\Omega}(\vth (\tau)-n)_+ \dx   \le \int_0^T\intom |g| \chi_{\{\vth\ge n\}} +\int_{\Omega}(\vth_0-n)_+ \dx.
    \end{aligned} 
\end{equation*}
Thus, using the fact that $\vth, g\in L^1(Q)$ and $\vth_0\in L^1(\Omega)$, the above inequality implies 
\begin{equation}\label{equire}
\lim_{n\to +\infty} \essup_{t\in (0,T)} \int_{\Omega} (\vth(t,x)-n)_+\dx=0,
\end{equation}
which implies the validity of \eqref{equi} since $\vth\chi_{\{\vth\geq n\}}/2\leq(\vth-n/2)\chi_{\{\vth\geq n/2\}}$ on $Q$ for $n>0$.

The chain of implications {\bf A)} $\implies$ {\bf E)} $\implies$
{\bf D)} $\implies$ {\bf A)} remains to be proved. 
We start by testing \eqref{iebR} with $\varphi=(1-(\cQ_m(\hat\vth)/\cQ_m(\vth))^\alpha)\chi_{(0,\tau)}$, where $\alpha>1/2$, $m>2\|\hat\vth\|_\infty$, $\tau\in [0,T]$ and the functions $\cQ_m$ and $\cH_{m,\alpha}$ are defined as follows
$$
\forall s>0: \cQ_m(s)=\int_0^s\kappa(\sigma)\cT'_{m,m/2}(\sigma)\dd \sigma,\quad\cH_{m,\alpha}(s)=\int_\mu^s\frac{\cT'_{m,m/2}(\sigma)}{\cQ_{m}(\sigma)^\alpha}\dd\sigma.
$$
Note that $\varphi\in L^\infty((0,T)\times\Omega)\cap L^2(0,T;W^{1,2}_0(\Omega))$, however, it is not as regular as required, particularly in time. By averaging in time, similarly as to get \eqref{eeB1.2}, we obtain for a.e. $\tau\in(0,T)$ and any $\alpha>1/2$, $\alpha\neq 1$, the following equality
\begin{equation}
\label{reem}
\begin{multlined}
\intom \cT_{m,m/2}(\vth(\tau))-\cQ_m(\hat\vth)^\alpha\cH_{m,\alpha}(\vth(\tau))-\big(\cT_{m,m/2}(\vth_0)-\cQ_m(\hat\vth)^\alpha\cH_{m,\alpha}(\vth_0)\big)\dx\\
-\int_0^\tau\intom \bu\cQ_m(\hat\vth)^\alpha\nabla\vth\frac{\cT'_{m,m/2}(\vth)}{\cQ_m(\vth)^\alpha}\dx\dt 
+\alpha\int_0^\tau\intom \cQ_m(\hat\vth)\left|\nabla \left(\frac{\cQ_m(\vth)}{\cQ_m(\hat\vth)}\right)\right|^2\left(\frac{\cQ_m(\vth)}{\cQ_m(\hat\vth)}\right)^{-1-\alpha}\dx\dt
\\
=\int_0^\tau\intom g\cT'_{m,m/2}(\vth) \left(1-\left(\frac{\cQ_m(\hat\vth)}{\cQ_m(\vth)}\right)^\alpha\right)\dx\dt+\int_0^\tau\intom (-\cT''_{m,m/2}(\vth)\kappa(\vth)|\nabla\vth|^2 \left(1-\left(\frac{\cQ_m(\hat\vth)}{\cQ_m(\vth)}\right)^\alpha\right)\dx\dt
\\
+\frac{\alpha}{\alpha-1}\int_0^\tau\intom\nabla\cQ_m(\hat\vth)\nabla\left( \left(\frac{\cQ_m(\hat\vth)}{\cQ_m(\vth)}\right)^{1-\alpha}-1\right)\dx\dt.
    \end{multlined}
\end{equation}
The calculations leading to the form of the elliptic term can be found, for instance, in \cite[Section~3.5]{Abbatiello2024}, while the remaining computation is straightforward. The last term in \eqref{reem} vanishes due to \eqref{hatcond} and the choice of sufficiently large $m$ so that $\cQ_m(\hat\vth)=\cQ(\hat\vth)$. In the case $\alpha=1$ the computation differs only in the form of the last term in \eqref{reem}; however, this term also vanishes. Hence, the argument applies to arbitrary $\alpha>1/2$.

Now we realize that due to \eqref{conditionstheta}, \eqref{vb2B}, \eqref{vb1}, and \eqref{vb2} we can pass to the limit as $m\to+\infty$ in the first and the second integral on the left-hand side of \eqref{reem} and in the first integral on the right-hand side, for a.e. $\tau\in(0,T)$, obtaining the corresponding terms appearing in \eqref{ree}.
The third integral on the right hand side of \eqref{reem} can be rewritten as 
\begin{equation}
\label{ree-ell}
\begin{multlined}
\alpha\int_0^\tau\intom \cQ_m(\hat\vth)\left|\nabla \left(\frac{\cQ_m(\vth)}{\cQ_m(\hat\vth)}\right)\right|^2 \left(\frac{\cQ_m(\vth)}{\cQ_m(\hat\vth)}\right)^{-1-\alpha}\dx\dt=
\alpha\int_0^\tau\intom \cQ(\hat\vth)\left|\nabla \left(\frac{\cQ(\vth)}{\cQ(\hat\vth)}\right )\right|^2 \left(\frac{\cQ(\vth)}{\cQ(\hat\vth)}\right)^{-1-\alpha}\dx\dt\\
+\alpha\int_0^\tau\intom \cQ_m(\hat\vth)\left|\nabla \left(\frac{\cQ_m(\vth)}{\cQ_m(\hat\vth)}\right)\right|^2 \left(\frac{\cQ_m(\vth)}{\cQ_m(\hat\vth)}\right)^{-1-\alpha}\chi_{\{\vth>m/2\}}\dx\dt\\
-
\alpha\int_0^\tau\intom \cQ(\hat\vth)\left|\nabla \left(\frac{\cQ(\vth)}{\cQ(\hat\vth)}\right) \right|^2 \left(\frac{\cQ(\vth)}{\cQ(\hat\vth)}\right)^{-1-\alpha}\chi_{\{\vth>m/2\}}\dx\dt.
\end{multlined}
\end{equation}
The last two integrals in \eqref{ree-ell} can be estimated by 
$$
C\int_0^T\intom \big(|\nabla\vth|^2\vth^{-1-\alpha}+|\nabla\hat\vth|^2(\vth^{-\alpha}+m^{-\alpha})\big)\chi_{\{\vth>m/2\}}
$$
which is a quantity that is independent of $\tau\in(0,T)$ and that converges to $0$ as $m\to+\infty$. Together, we can rewrite \eqref{reem} for a.e. $\tau\in(0,T)$ as 
\begin{equation}
    \label{reem1}
\begin{multlined}
\intom \vth(\tau)-\cQ(\hat\vth)^\alpha\cH_\alpha(\vth(\tau))-\big(\vth_0-\cQ(\hat\vth)^\alpha\cH_\alpha(\vth_0)\big)\dx-\int_0^\tau\intom \bu\cQ(\hat\vth)^\alpha\nabla\cH_\alpha(\vth)\dx\dt \\
+\alpha\int_0^\tau\intom \cQ(\hat\vth)\left|\nabla \left(\frac{\cQ(\vth)}{\cQ(\hat\vth)}\right)\right|^2 \left(\frac{\cQ(\vth)}{\cQ(\hat\vth)}\right)^{-1-\alpha}\dx\dt-\int_0^\tau\intom g \left(1-\left(\frac{\cQ(\hat\vth)}{\cQ(\vth)}\right)^\alpha\right)\dx\dt\\
= \lim_{m\to\infty}
\int_0^\tau\intom (-\cT''_{m,m/2}(\vth)\kappa(\vth)|\nabla\vth|^2 \left(1-\left(\frac{\cQ_m(\hat\vth)}{\cQ_m(\vth)}\right)^\alpha\right)\dx\dt+\lim_{m\to\infty}I_m(\tau),
\end{multlined}
\end{equation}
where $I_m(\tau)\to 0$ as $m\to+\infty$ for a.e. $\tau\in (0,T)$.

To prove {\bf A)} $\implies$ {\bf E)} we estimate
$$
\int_0^\tau\intom (-\cT''_{m,m/2}(\vth)\kappa(\vth)|\nabla\vth|^2 \left(1-\left(\frac{\cQ_m(\hat\vth)}{\cQ_m(\vth)}\right)^\alpha\right)
\leq
C\int_0^T\intom \frac{|\nabla\vth|^2}{m}\chi_{\{\vth\in(m/2,3m/2)\}}
$$
by \eqref{propertiestke} and the estimate $\vth\geq\vth_0$. It remains to use \eqref{Equi1} for limit passage as $m\to+\infty$ in \eqref{reem1} to get the relative energy equality in {\bf E)}. 

The implication {\bf E)} $\implies$ {\bf D)} is obvious. 

It remains to show {\bf D)} $\implies$ {\bf A)}.
Using \eqref{propertiestke}, we estimate
\[
\begin{multlined}
\int_0^\tau\intom (-\cT''_{m,m/2}(\vth)\kappa(\vth)|\nabla\vth|^2 \left(1-\left(\frac{\cQ_m(\hat\vth)}{\cQ_m(\vth)}\right)^\alpha\right)
\geq
C\int_0^\tau\intom \frac{|\nabla\vth|^2}{m}\chi_{\{\vth\in(3m/4,5m/4)\}}(1-c\frac{\cQ(\|\hat\vth\|_\infty)}m)
\\
\geq C\int_0^\tau\intom \frac{|\nabla\vth|^2}{m}\chi_{\{\vth\in(3m/4,5m/4)\}}
\end{multlined}\]
for sufficiently large $m$. The left-hand side of~\eqref{reem1} is non-positive for some $\alpha>1/2$ due to the relative energy inequality in~{\bf D)}. Passing to the limit in~\eqref{reem1} and using the estimate above, we have for a.e. $\tau\in(0,T)$
$$
\lim_{m\to+\infty}\int_0^\tau\intom \frac{|\nabla\vth|^2}{m}\chi_{\{\vth\in(3m/4,5m/4)\}}=0.
$$
This concludes the proof of Theorem~\ref{maintheoremB}.


\end{document}